\documentclass[12pt,reqno]{amsart}

\newcommand{\df}{\dfrac}
\newcommand{\tf}{\tfrac}

\renewcommand{\(}{\left\(}
\renewcommand{\)}{\right\)}
\renewcommand{\[}{\left\[}
\renewcommand{\]}{\right\]}
\let\dotlessi=\i
\renewcommand{\i}{\infty}

\numberwithin{equation}{section}
 \theoremstyle{plain}
\newtheorem{theorem}{Theorem}[section]
\newtheorem{lemma}[theorem]{Lemma}
\newtheorem{corollary}[theorem]{Corollary}

   \makeatletter
\def\proof{\@ifnextchar[{\@oproof}{\@nproof}}
\def\@oproof[#1][#2]{\trivlist\item[\hskip\labelsep\textit{#2 Proof of\
#1.}~]\ignorespaces}
\def\@nproof{\trivlist\item[\hskip\labelsep\textit{Proof.}~]\ignorespaces}

\makeatother

\begin{document}
\title[Character analogues]{Character analogues of Ramanujan type integrals involving the Riemann $\Xi$-function}
\author{Atul Dixit}
\address{Department of Mathematics, University of Illinois, 1409 West Green
Street, Urbana, IL 61801, USA} \email{aadixit2@illinois.edu}
\maketitle
\begin{abstract}
A new class of integrals involving the product of $\Xi$-functions associated with primitive Dirichlet characters is considered. These integrals give rise to transformation formulas of the type $F(z, \alpha,\chi)=F(-z, \beta,\overline{\chi})=F(-z,\alpha,\overline{\chi})=F(z,\beta,\chi)$, where $\alpha\beta=1$. New character analogues of transformation formulas of Guinand and Koshliakov as well as those of a formula of Ramanujan and its recent generalization are shown as particular examples. Finally, character analogues of a conjecture of Ramanujan, Hardy and Littlewood involving infinite series of M\"{o}bius functions are derived.
\end{abstract}
\section{Introduction}
Modular transformations are ubiquitous in Ramanujan's Notebooks \cite{nb} and his Lost Notebook \cite{lnb}. Ramanujan usually expressed them in a symmetric way and they were valid under the conditions $\alpha\beta=\pi$ or $\alpha\beta=\pi^2$ etc. In the same spirit, on page 220 in one of the manuscripts of S.~Ramanujan in the handwriting of G.N.~Watson published in \cite{lnb}, one finds the following beautiful claim.
\begin{theorem}\label{entry1} Define
\begin{equation*}\label{w1.27}
\lambda(x):=\psi(x)+\df{1}{2x}-\log x,
\end{equation*}
where
\begin{equation*}\label{w1.15b}
\psi(x):=\df{\Gamma^\prime(x)}{\Gamma(x)}=-\gamma-\sum_{m=0}^{\infty}\left(\df{1}{m+x}-\df{1}{m+1}\right)
\end{equation*}
is the logarithmic derivative of the Gamma function.
Let the Riemann $\xi$-function be defined by
\begin{equation*}
\xi(s):=(s-1)\pi^{-\tf{1}{2}s}\Gamma(1+\tf{1}{2}s)\zeta(s),
\end{equation*}
and let 
\begin{equation*}\label{xif}
\Xi(t):=\xi(\tf{1}{2}+it)
\end{equation*}
be the Riemann $\Xi$-function.
If $\alpha$ and $\beta$ are positive numbers such that $\alpha\beta=1$, then
\begin{align}\label{w1.26}
\sqrt{\alpha}\bigg(\frac{\gamma-\log(2\pi\alpha)}{2\alpha}&+\sum_{k=1}^{\infty}\lambda(k\alpha)\bigg)
=\sqrt{\beta}\bigg(\frac{\gamma-\log(2\pi\beta)}{2\beta}+\sum_{k=1}^{\infty}\lambda(k\beta)\bigg)\nonumber\\
&=-\frac{1}{\pi^{3/2}}\int_0^{\infty}\left|\Xi\left(\frac{1}{2}t\right)\Gamma\left(\frac{-1+it}{4}\right)\right|^2
\frac{\cos\left(\tfrac{1}{2}t\log\alpha\right)}{1+t^2}\, dt,
\end{align}
where $\gamma$ denotes Euler's constant.
\end{theorem}
This identity is of a special kind since not only does it contain a modular transformation, but also a beautiful integral involving the Riemann $\Xi$-function. In fact, the invariance of the integral in (\ref{w1.26}) under the map $\alpha\to\beta$ establishes the equality of the first and the second expressions in (\ref{w1.26}).
This idea is used in \cite{bcbad} to prove the above claim of Ramanujan and later in \cite{dixit, series, transf} to obtain many transformation formulas of the type $F(\alpha)=F(\beta)$ or $F(z,\alpha)=F(z,\beta)$, where $\alpha\beta=1$, where an integral involving the Riemann $\Xi$-function is always linked to them. This then gives new identities involving infinite series of Hurwitz zeta function as well as extensions of some well-known formulas of A.P.~Guinand \cite{guinand} and N.S.~Koshliakov \cite{koshliakov} (see \cite{transf} for details). For example, we mention the following generalization of Theorem \ref{entry1} found in \cite{dixit, transf}. 
\begin{theorem}\label{sechur}
Let $-1<$ \textup{Re} $z<1$. Define $\varphi(z,x)$ by
\begin{equation*}\label{dvarphi}
\varphi(z,x)=\zeta(z+1,x)-\frac{x^{-z}}{z}-\frac{1}{2}x^{-z-1},
\end{equation*}
where $\zeta(z,x)$ denotes the Hurwitz zeta function. Then if $\alpha$ and $\beta$ are any positive numbers such that $\alpha\beta=1$,
\begin{align}\label{mainneq2}
&\alpha^{\frac{z+1}{2}}\left(\sum_{n=1}^{\infty}\varphi(z,n\alpha)-\frac{\zeta(z+1)}{2\alpha^{z+1}}-\frac{\zeta(z)}{\alpha z}\right)=\beta^{\frac{z+1}{2}}\left(\sum_{n=1}^{\infty}\varphi(z,n\beta)-\frac{\zeta(z+1)}{2\beta^{z+1}}-\frac{\zeta(z)}{\beta z}\right)\nonumber\\
&=\frac{8(4\pi)^{\frac{z-3}{2}}}{\Gamma(z+1)}\int_{0}^{\infty}\Gamma\left(\frac{z-1+it}{4}\right)\Gamma\left(\frac{z-1-it}{4}\right)
\Xi\left(\frac{t+iz}{2}\right)\Xi\left(\frac{t-iz}{2}\right)\frac{\cos\left( \tf{1}{2}t\log\alpha\right)}{(z+1)^2+t^2}\, dt,
\end{align}
where $\Xi(t)$ is the Riemann $\Xi$-function. 
\end{theorem}
Another example of a transformation formula of the type $F(z,\alpha)=F(z,\beta)$ is the extended version of Guinand's formula \cite[Theorem 1.4]{transf} given below.
\begin{theorem}\label{mainnn}
Let $K_{\nu}(s)$ denote the modified Bessel function of order $\nu$, let $\gamma$ denote Euler's constant and let $\sigma_k(n)=\sum_{d|n}d^k$. Let $-1<$ \textup{Re} $z<1$. Then if $\alpha$ and $\beta$ are positive numbers such that $\alpha\beta=1$, we have
\begin{align}\label{mainagain3}
&\sqrt{\alpha}\left(\alpha^{\frac{z}{2}-1}\pi^{\frac{-z}{2}}\Gamma\left(\frac{z}{2}\right)\zeta(z)+\alpha^{-\frac{z}{2}-1}\pi^{\frac{z}{2}}\Gamma\left(\frac{-z}{2}\right)\zeta(-z)-4\sum_{n=1}^{\infty}\sigma_{-z}(n)n^{z/2}K_{\frac{z}{2}}\left(2n\pi\alpha\right)\right)\nonumber\\
&=\sqrt{\beta}\left(\beta^{\frac{z}{2}-1}\pi^{\frac{-z}{2}}\Gamma\left(\frac{z}{2}\right)\zeta(z)+\beta^{-\frac{z}{2}-1}\pi^{\frac{z}{2}}\Gamma\left(\frac{-z}{2}\right)\zeta(-z)-4\sum_{n=1}^{\infty}\sigma_{-z}(n)n^{z/2}K_{\frac{z}{2}}\left(2n\pi\beta\right)\right)\nonumber\\
&=-\frac{32}{\pi}\int_{0}^{\infty}\Xi\left(\frac{t+iz}{2}\right)\Xi\left(\frac{t-iz}{2}\right)\frac{\cos\left(\frac{1}{2}t\log\alpha\right)}{(t^2+(z+1)^2)(t^2+(z-1)^2)}\, dt.
\end{align}
\end{theorem}
Letting $z\to 0$ in (\ref{mainagain3}) then gives the following extended version of Koshliakov's formula \cite{series}.
\begin{theorem}\label{koshfe2r}
Let $d(n)$ denote the number of positive divisors of $n$ and let $K_{0}(n)$ denote the modified Bessel function of order $0$. If $\alpha$ and $\beta$ are positive numbers such that $\alpha\beta=1$, then
\begin{align*}\label{koshfe}
&\sqrt{\alpha}\left(\frac{\gamma-\log (4\pi\alpha)}{\alpha}-4\sum_{n=1}^{\infty}d(n)K_{0}(2\pi n\alpha)\right) =\sqrt{\beta}\left(\frac{\gamma-\log (4\pi\beta)}{\beta}-4\sum_{n=1}^{\infty}d(n)K_{0}(2\pi n\beta)\right)\nonumber\\
&=-\frac{32}{\pi}\int_{0}^{\infty}\frac{\left(\Xi\left(\frac{t}{2}\right)\right)^2\cos\left(\frac{1}{2}t\log\alpha\right)\, dt}{(1+t^2)^2}.
\end{align*}
\end{theorem}
By an `extended version', we mean that the original identity known before is linked to an integral involving the Riemann $\Xi$-function. It must be mentioned here that the formulas of Guinand and Koshliakov were discovered earlier by Ramanujan (see \cite{bls}), and are present in the Lost Notebook \cite[p.~253--254]{lnb}. 

After Ramanujan, Koshliakov was another person to do significant research in this area. Apart from using contour integration, Mellin transforms and several summation formulas that he developed, he frequently used a method similar to that developed by Ramanujan in \cite{riemann}, \cite[pp.~72--77]{cp} to obtain old and new transformation formulas of the form $F(\alpha)=F(\beta)$, where $\alpha\beta=k$ for some constant $k$. He also obtained deep generalizations of many well-known formulas of Ramanujan and of G.H.~Hardy \cite[Equation (2)]{ghh}, some of them being analogues in rational and number fields. See \cite{koshliakov6, koshliakov3, koshliakov1, koshliakov7, koshqua}. In \cite{koshliakov5, koshxi}, he used Fourier's integral theorem to obtain expressions for the Riemann $\Xi$-function, a method also enunciated by Ramanujan \cite{riemann}. Around the same time, W.L.~Ferrar \cite{ferrar} also worked on transformation formulas of the above kind.

As can be seen from (\ref{w1.26}), the general form of the integrals giving rise to formulas of the type $F(\alpha)=F(\beta)$ where $\alpha\beta=1$ is 
\begin{equation*}\label{tint}
\int_{0}^{\infty}f\left(\frac{t}{2}\right)\Xi\left(\frac{t}{2}\right)\cos \mu t\, dt, 
\end{equation*}
for $\mu$ real and $f(t)=\phi(it)\phi(-it)$, where $\phi$ is analytic in $t$ as a function of a real variable. This integral is mentioned in \cite[p.~35]{titch}.
Similarly, from (\ref{mainneq2}) and (\ref{mainagain3}), it is clear that the general form of the integrals giving rise to identities of the type $F(z,\alpha)=F(z,\beta)$ where $\alpha\beta=1$ is
\begin{equation}\label{genint}
\int_{0}^{\infty}f\left(z,\frac{t}{2}\right)\Xi\left(\frac{t+iz}{2}\right)\Xi\left(\frac{t-iz}{2}\right)\cos \mu t\, dt,
\end{equation}
for $\mu$ real and $f(z,t)=\phi(z,it)\phi(z,-it)$, where $\phi$ is analytic in $t$ as a function of a real variable as well as analytic in $z$ in some complex domain. An integral of this kind was first introduced by Ramanujan \cite{riemann}. 

In this article, we find character analogues of all of the above-mentioned theorems. The character analogue of Guinand's, and hence Koshliakov's formula given here differs from the ones established in \cite{bds}. Throughout this article, we will be concerned with the principal branch of the logarithm. Also, we work only with a primitive, non-principal Dirichlet character $\chi$ modulo $q$. It is easy to see that its conjugate character $\overline{\chi}$ is also a primitive, non-principal character modulo $q$ and $\overline{\chi}$ is even (odd) if and only if $\chi$ is even (resp. odd). Let $L(s,\chi)$ denote the Dirichlet $L$-function defined by
$L(s,\chi)=\sum_{n=1}^{\infty}\chi(n)/n^s$ for Re $s>1$. This series converges conditionally for $0<$ Re $s<1$. Also, it can be analytically continued to an entire function of $s$. Let $G(\chi):=G(1,\chi)$, where $G(n,\chi)$ is the Gauss sum defined by
\begin{equation*}\label{gauss}
G(n,\chi):=\sum_{m=1}^{q}\chi(m)e^{2\pi imn/q}.
\end{equation*}
We know that \cite[p.~168]{apostol}
\begin{equation}\label{sgs}
|G(\chi)|^{2}=q
\end{equation} 
and it is easy to see that
\begin{equation}\label{eo}
\overline{G(\chi)}=
\begin{cases}
G(\overline{\chi}), \quad\mbox{for}\hspace{1mm}\chi\hspace{1mm}\mbox{even},\\
-G(\overline{\chi}), \quad\mbox{for}\hspace{1mm}\chi\hspace{1mm}\mbox{odd}.\\
\end{cases}
\end{equation}
Define $b$ as follows:
\begin{equation}\label{a}
b=
\begin{cases}
0, \quad\chi(-1)=1,\\
1, \quad\chi(-1)=-1.\\
\end{cases}
\end{equation}
Then the function $\xi(s,\chi)$ is defined by
\begin{equation}\label{xic}
\xi(s,\chi):=\left(\frac{\pi}{q}\right)^{-(s+b)/2}\Gamma\left(\frac{s+b}{2}\right)L(s,\chi),
\end{equation}
and the analogue of the Riemann $\Xi$-function for Dirichlet characters is then defined as
\begin{equation}\label{Xic}
\Xi(t,\chi):=\xi\left(\frac{1}{2}+it,\chi\right).
\end{equation}
$L$-functions satisfy the functional equation \cite[p.~263]{apostol}
\begin{equation}\label{funl}
L(1-s,\chi)=\frac{q^{s-1}\Gamma(s)}{(2\pi)^{s}}\left(e^{-\pi is/2}+\chi(-1)e^{\pi is/2}\right)G(\chi)L(s,\overline{\chi}),
\end{equation}
which can be rephrased in terms of $\xi(s,\chi)$ as \cite{dav}
\begin{equation}\label{fund}
\xi(1-s,\overline{\chi})=\epsilon (\chi)\xi(s,\chi),
\end{equation}
where $\epsilon (\chi)=i^{b}q^{1/2}/G(\chi)$. By (\ref{sgs}), $|\epsilon (\chi)|=1$. Next, we note Stirling's formula in a vertical strip $\alpha\leq\sigma\leq\beta$, $s=\sigma+it$, namely,
\begin{equation}\label{strivert}
|\Gamma(s)|=(2\pi)^{\tf{1}{2}}|t|^{\sigma-\tf{1}{2}}e^{-\tf{1}{2}\pi |t|}\left(1+O\left(\frac{1}{|t|}\right)\right)
\end{equation}
uniformly as $|t|\to\infty$. Now using (\ref{funl}) and the fact \cite[p.~82]{dav} that $|L(s,\chi)|=O(q|t|)$ for Re $s\geq 1/2$, we can easily see that for Re $s\geq-\delta$, $\delta>0$, we have
\begin{equation}\label{lbound}
L(s,\chi)=O\left(q^{\frac{3}{2}+\delta}|t|^{\frac{3}{2}+\delta}\right).
\end{equation}
We will subsequently use this result. 

Transformation formulas involving Dirichlet characters of the form
\begin{equation*}
\sum_{n=1}^{\infty}\chi(n)f(n)=\sum_{n=1}^{\infty}\overline{\chi}(n)g(n),
\end{equation*}
where
\begin{equation*}
g(x)=
\begin{cases}
\displaystyle\frac{2G(\chi)}{q}\int_{0}^{\to\infty}\cos\left(\frac{2\pi xt}{q}\right)f(t)\, dt,\quad\mbox{for}\hspace{1mm}\chi(-1)=1,\\[2ex]
\displaystyle\frac{-2iG(\chi)}{q}\int_{0}^{\to\infty}\sin\left(\frac{2\pi xt}{q}\right)f(t)\, dt,\quad\mbox{for}\hspace{1mm}\chi(-1)=-1,
\end{cases}
\end{equation*}
were considered by Guinand \cite[Theorems 4--5]{guipoi}, though he did not give any particular examples. Here, we derive a character analogue of the integral in (\ref{genint}). Its general form is
\begin{equation}\label{gfica}
\int_{0}^{\infty}f\left(z,\frac{t}{2}\right)\Xi\left(\frac{t+iz}{2},\overline{\chi}\right)\Xi\left(\frac{t-iz}{2},\chi\right)\cos\left(\frac{1}{2}t\log\alpha\right)\, dt,
\end{equation}
where $f$ is an even function of both the variables $z$ and $t$. These integrals give rise to transformation formulas of the type $F(z,\alpha,\chi)=F(-z,\beta,\overline{\chi})=F(-z,\alpha,\overline{\chi})=F(z,\beta,\chi)$. Then via Fourier's integral theorem, one may be able to obtain integral representations for $\Xi\left((t+iz)/2,\overline{\chi}\right)\Xi\left((t-iz)/2,\chi\right)$ which are of independent interest.
The character analogue of Theorem \ref{mainnn} is as follows.
\begin{theorem}\label{koshfec}
Let $-1<$ Re $z<1$ and let $\chi$ denote a primitive, non-principal character modulo $q$. Let the number $b$ be defined as in \textup{(\ref{a})}. Let $K_{\nu}(z), d(n)$ and $\gamma$ be defined as before and let $\alpha$ and $\beta$ be positive numbers such that $\alpha\beta=1$. If
\begin{equation*}\label{fdefe}
F(z,\alpha,\chi):=\alpha^{b+\frac{1}{2}}\sum_{n=1}^{\infty}\chi(n)n^{-\frac{z}{2}+b}\bigg(\sum_{d|n}\overline{\chi}^{2}(d)d^{z}\bigg)K_{-\frac{z}{2}}\left(\frac{2\pi n\alpha}{q}\right),
\end{equation*}
then
\begin{align}\label{koshfec1}
F(z,\alpha,\chi)&=F(-z,\beta,\overline{\chi})=F(-z,\alpha,\overline{\chi})=F(z,\beta,\chi)\nonumber\\
&=\frac{1}{8\pi}\int_{0}^{\infty}\Xi\left(\frac{t+iz}{2},\overline{\chi}\right)\Xi\left(\frac{t-iz}{2},\chi\right)\cos\left(\frac{1}{2}t\log\alpha\right)\, dt.
\end{align}
\end{theorem}
Define $\psi(a,\chi)$ by
\begin{equation}\label{charpsi}
\psi(a,\chi)=-\sum_{n=1}^{\infty}\frac{\chi(n)}{n+a},
\end{equation}
where $a\in\mathbb{C}$ is a non-negative integer. For a real character $\chi$, this agrees with the character analogue of the psi function obtained by the logarithmic differentiation of the following Weierstrass product form of the character analogue of the gamma function for real characters derived by Berndt \cite{berndt}:
\begin{equation*}\label{chargw}
\Gamma(a,\chi)=e^{-aL(1,\chi)}\prod_{n=1}^{\infty}\left(1+\frac{a}{n}\right)^{-\chi(n)}e^{a\chi(n)/n}.
\end{equation*}
The character analogue of the Hurwitz zeta function $\zeta(z,a)$ is given by \cite[Ex. 3.2]{bcberndt}
\begin{equation}\label{charhur}
L(z,a,\chi)=\sum_{n=1}^{\infty}\frac{\chi(n)}{(n+a)^{z}},
\end{equation}
valid for Re $z>0$, and provided $a\in\mathbb{C}$ is a non-negative integer.
The above character analogue of the Hurwitz zeta function can also be obtained as the special case when $x=0$ of the function $L(z,x,a,\chi)$ defined in \cite{berndt} by 
\begin{equation*}\label{berndtg}
L(z,x,a,\chi):={\sum_{n=0}^{\infty}}'e^{2\pi inx/k}\chi(n)(n+a)^{-z},
\end{equation*}
where the prime indicates that the term corresponding to $n=-a$ is omitted if $a$ is a negative integer and $\chi(a)\neq 0$. As shown in \cite{berndt}, $L(z,x,a,\chi)$ converges for Re $z>0$ if $x$ is not an integer, or if $x$ is an integer and gcd$(x,k)>1$. If $x$ is an integer and gcd$(x,k)=1$, the series converges for Re $z>1$. For mean value properties of $L(z,a,\chi)$ and asymptotic formulas, see the recent paper \cite{myx}.
The character analogues of Theorem \ref{sechur} are given below.
\begin{theorem}\label{ramgenec}
Let $\chi$ denote an even, primitive, non-principal character modulo $q$. Let $-1<$ Re $z<1$ and let $L(z,a,\chi)$ be defined as in \textup{(\ref{charhur})}. Define $T(z,\alpha,\chi)$ by
\begin{equation}\label{te}
T(z,\alpha,\chi):=\frac{\alpha^{z/2}q^{z/2}\Gamma(z+1)}{2^{z}\pi^{z/2}G(\chi)},
\end{equation}
and $\Omega(z,t)$ by
\begin{align}\label{omdef}
\Omega(z,t)&:=((z+1)^2+t^2)\Gamma\bigg(\frac{-z-1+it}{4}\bigg)\Gamma\bigg(\frac{-z-1-it}{4}\bigg)\nonumber\\
&\quad+((z-1)^2+t^2)\Gamma\bigg(\frac{z-1+it}{4}\bigg)\Gamma\bigg(\frac{z-1-it}{4}\bigg).
\end{align}
Then if $\alpha$ and $\beta$ are positive numbers such that $\alpha\beta=1$,
{\allowdisplaybreaks\begin{align}\label{ramgenec1}
&\sqrt{\alpha}\bigg[T(z,\alpha,\chi)\sum_{n=1}^{\infty}\chi(n)L\bigg(z+1,n\alpha,\chi\bigg)+T(-z,\alpha,\overline{\chi})\sum_{n=1}^{\infty}\overline{\chi}(n)L\bigg(-z+1,n\alpha,\overline{\chi}\bigg)\bigg]\nonumber\\
&=\sqrt{\beta}\bigg[T(-z,\beta,\overline{\chi})\sum_{n=1}^{\infty}\overline{\chi}(n)L\bigg(-z+1,n\beta ,\overline{\chi}\bigg)+T(z,\beta,\chi)\sum_{n=1}^{\infty}\chi(n)L\bigg(z+1,n\beta,\chi\bigg)\bigg]\nonumber\\
&=\frac{1}{64\pi^{3/2}q}\int_{0}^{\infty}\Omega(z,t)\Xi\left(\frac{t+iz}{2},\overline{\chi}\right)\Xi\left(\frac{t-iz}{2},\chi\right)\cos\left(\frac{1}{2}t\log\alpha\right)\, dt.
\end{align}}
\end{theorem}
\begin{theorem}\label{ramgenoc}
Let $\chi$ denote an odd, primitive, non-principal character modulo $q$. Let $-1<$ Re $z<1$ and let $L(z,a,\chi)$ be defined as in \textup{(\ref{charhur})} and let $T(z,\alpha,\chi)$ be defined as in \textup{(\ref{te})}. Define $\Lambda(z,t)$ by
\begin{equation}\label{omdefo}
\Lambda(z,t):=\Gamma\bigg(\frac{z+1+it}{4}\bigg)\Gamma\bigg(\frac{z+1-it}{4}\bigg)+\Gamma\bigg(\frac{-z+1+it}{4}\bigg)\Gamma\bigg(\frac{-z+1-it}{4}\bigg).
\end{equation}
Then if $\alpha$ and $\beta$ are positive numbers such that $\alpha\beta=1$,
\begin{align}\label{ramgenoc1}
&\sqrt{\alpha}\bigg[T(z,\alpha,\chi)\sum_{n=1}^{\infty}\chi(n)L\bigg(z+1,n\alpha,\chi\bigg)+T(-z,\alpha,\overline{\chi})\sum_{n=1}^{\infty}\overline{\chi}(n)L\bigg(-z+1,n\alpha,\overline{\chi}\bigg)\bigg]\nonumber\\
&=\sqrt{\beta}\bigg[T(-z,\beta,\overline{\chi})\sum_{n=1}^{\infty}\overline{\chi}(n)L\bigg(-z+1,n\beta ,\overline{\chi}\bigg)+T(z,\beta,\chi)\sum_{n=1}^{\infty}\chi(n)L\bigg(z+1,n\beta,\chi\bigg)\bigg]\nonumber\\
&=\frac{1}{4\pi^{1/2}iq^{2}}\int_{0}^{\infty}\Lambda(z,t)\Xi\left(\frac{t+iz}{2},\overline{\chi}\right)\Xi\left(\frac{t-iz}{2},\chi\right)\cos\left(\frac{1}{2}t\log\alpha\right)\, dt.
\end{align}
\end{theorem}
In \cite[p.~156, Section 2.5]{hl}, Hardy and Littlewood discuss the following interesting identity suggested to them by work of Ramanujan. 
\begin{theorem}\label{rhl}
Let $\mu(n)$ denote the M\"obius function. Let $\alpha$ and $\beta$ be two positive numbers such that $\alpha\beta=1$. Assume that the series $\sum_{\rho}\left(\Gamma{\left(\frac{1-\rho}{2}\right)}/\zeta^{'}(\rho)\right)a^{\rho}$ converges, where $\rho$ denotes a non-trivial zero of the Riemann zeta function and $a$ denotes a positive real number, and that the non-trivial zeros of $\zeta(s)$ are simple. Then
\begin{align}\label{mr}
&\sqrt{\alpha}\sum_{n=1}^{\infty}\frac{\mu(n)}{n}e^{-\pi\alpha^2/n^2}-\frac{1}{4\sqrt{\pi}\sqrt{\alpha}}\sum_{\rho}\frac{\Gamma{\left(\frac{1-\rho}{2}\right)}}{\zeta^{'}(\rho)}\pi^{\frac{\rho}{2}}\alpha^{\rho}\nonumber\\
&=\sqrt{\beta}\sum_{n=1}^{\infty}\frac{\mu(n)}{n}e^{-\pi\beta^2/n^2}-\frac{1}{4\sqrt{\pi}\sqrt{\beta}}\sum_{\rho}\frac{\Gamma{\left(\frac{1-\rho}{2}\right)}}{\zeta^{'}(\rho)}\pi^{\frac{\rho}{2}}\beta^{\rho}.
\end{align}
\end{theorem}
The original formulation of the above identity is slightly different in \cite{hl} but can readily be seen to be equivalent to (\ref{mr}). See also \cite[p.~470]{berndt1}, \cite[p.~143]{kp} and \cite[p.~219, Section 9.8]{titch} for discussions on this identity. Based on certain assumptions, the character analogues of (\ref{mr}) for even and odd primitive Dirichlet characters, which furnish two examples of transformation formulas of the form $F(\alpha,\chi)=F(\beta,\overline{\chi})$, are derived here and are as follows.
\begin{theorem}\label{rhlao}
Let $\chi$ be an odd, primitive character modulo $q$, and let $\alpha$ and $\beta$ be two positive numbers such that $\alpha\beta=1$. Assume that the series $\sum_{\rho}\frac{\pi^{\rho/2}\alpha^{\rho}\Gamma{((2-\rho)/2)}}{q^{\rho/2}L{'}(\rho,\chi)}$ and $\sum_{\rho}\frac{\pi^{\rho/2}\beta^{\rho}\Gamma{((2-\rho)/2)}}{q^{\rho/2}L{'}(\rho,\overline{\chi})}$ converge, where $\rho$ denotes a non-trivial zero of $L(s,\chi)$ and $L(s,\overline{\chi})$ respectively, and that the non-trivial zeros of the associated Dirichlet L-functions are simple. Then
\begin{align}\label{mrao}
&\alpha\sqrt{\alpha}\sqrt{G(\chi)}\left(\sum_{n=1}^{\infty}\frac{\chi(n)\mu(n)}{n^2}e^{-\frac{\pi\alpha^2}{qn^2}}-\frac{q}{4\pi\alpha^2}\sum_{\rho}\frac{\Gamma{\left(\frac{2-\rho}{2}\right)}}{L{'}(\rho,\chi)}\left(\frac{\pi}{q}\right)^{\frac{\rho}{2}}\alpha^{\rho}\right)\nonumber\\
&=\beta\sqrt{\beta}\sqrt{G(\overline{\chi})}\left(\sum_{n=1}^{\infty}\frac{\overline{\chi}(n)\mu(n)}{n^2}e^{-\frac{\pi\beta^2}{qn^2}}-\frac{q}{4\pi\beta^2}\sum_{\rho}\frac{\Gamma{\left(\frac{2-\rho}{2}\right)}}{L{'}(\rho,\overline{\chi})}\left(\frac{\pi}{q}\right)^{\frac{\rho}{2}}\beta^{\rho}\right).
\end{align}
\end{theorem}
\begin{theorem}\label{rhlae}
Let $\chi$ be an even, primitive character modulo $q$, and let $\alpha$ and $\beta$ be two positive numbers such that $\alpha\beta=1$. Assume that the series $\sum_{\rho}\frac{\pi^{\rho/2}\alpha^{\rho}\Gamma{((2-\rho)/2)}}{q^{\rho/2}L{'}(\rho,\chi)}$ and $\sum_{\rho}\frac{\pi^{\rho/2}\beta^{\rho}\Gamma{((2-\rho)/2)}}{q^{\rho/2}L{'}(\rho,\overline{\chi})}$ converge, where $\rho$ denotes a non-trivial zero of $L(s,\chi)$ and $L(s,\overline{\chi})$ respectively, and that the non-trivial zeros of the associated Dirichlet L-functions are simple. Then
\begin{align}\label{mrae}
&\sqrt{\alpha}\sqrt{G(\chi)}\left(\sum_{n=1}^{\infty}\frac{\chi(n)\mu(n)}{n}e^{-\frac{\pi\alpha^2}{qn^2}}-\frac{\sqrt{q}}{4\sqrt{\pi}\alpha}\sum_{\rho}\frac{\Gamma{\left(\frac{1-\rho}{2}\right)}}{L{'}(\rho,\chi)}\left(\frac{\pi}{q}\right)^{\frac{\rho}{2}}\alpha^{\rho}\right)\nonumber\\
&=\sqrt{\beta}\sqrt{G(\overline{\chi})}\left(\sum_{n=1}^{\infty}\frac{\overline{\chi}(n)\mu(n)}{n}e^{-\frac{\pi\beta^2}{qn^2}}-\frac{\sqrt{q}}{4\sqrt{\pi}\beta}\sum_{\rho}\frac{\Gamma{\left(\frac{1-\rho}{2}\right)}}{L{'}(\rho,\overline{\chi})}\left(\frac{\pi}{q}\right)^{\frac{\rho}{2}}\beta^{\rho}\right).
\end{align}
\end{theorem}
This paper is organized as follows. In Section 2, we give a complex integral representation of (\ref{gfica}) that is used in subsequent sections. In Section 3, we prove Theorem \ref{koshfec}. Then in Section 4, we compute the inverse Mellin transforms and asymptotic expansions of certain functions which are subsequently used in Section 5. Section 5 is devoted to proofs of Theorems \ref{ramgenec} and \ref{ramgenoc}. Character analogues of Ramanujan's transformation formula (Theorem \ref{entry1}) are derived as special cases of these theorems. We conclude this section with some curious results on certain double series being always real. In Section 6, we present proofs of Theorems \ref{rhlao} and \ref{rhlae}. Finally we conclude with some open problems in Section 7.
\section{A complex integral representation of (\ref{gfica})}
In this section, we give a formal way of transforming an integral involving a character analogue of Riemann's $\Xi$-function into an equivalent complex integral which allows us to use residue calculus and Mellin transform techniques for its evaluation.
\begin{theorem}
Let 
\begin{equation}\label{finphi}
f(z,t)=\frac{\phi(z,it)\phi(z,-it)+\phi(-z,it)\phi(-z,-it)}{2},
\end{equation} 
where $\phi$ is analytic in $t$ as a function of a real variable and analytic in $z$ in some complex domain. Let $y=e^{\mu}$ with $\mu$ real. Then, under the assumption that the integral on the left side below converges,
\begin{align}\label{genint2}
&\int_{0}^{\infty}f(z,t)\Xi\left(t+\frac{iz}{2},\overline{\chi}\right)\Xi\left(t-\frac{iz}{2},\chi\right)\cos\mu t\, dt\nonumber\\
&=\frac{1}{4i\sqrt{y}}\int_{\frac{1}{2}-i\infty}^{\frac{1}{2}+i\infty}\bigg(\phi\left(z,s-\frac{1}{2}\right)\phi\left(z,\frac{1}{2}-s\right)+\phi\left(-z,s-\frac{1}{2}\right)\phi\left(-z,\frac{1}{2}-s\right)\bigg)\nonumber\\
&\quad\quad\quad\quad\quad\quad\quad\times\xi\left(s-\frac{z}{2},\overline{\chi}\right)\xi\left(s+\frac{z}{2},\chi\right)y^{s}\, ds.
\end{align}
\end{theorem}
\begin{proof}
Let 
\begin{equation*}\label{icnn}
I(z,\mu,\chi):=\int_{0}^{\infty}f(z,t)\Xi\left(t+\frac{iz}{2},\overline{\chi}\right)\Xi\left(t-\frac{iz}{2},\chi\right)\cos\mu t\, dt.
\end{equation*}
Then
\begin{align}\label{ic}
I(z,\mu,\chi)&=\frac{1}{2}\bigg(\int_{0}^{\infty}f(z,t)\Xi\left(t+\frac{iz}{2},\overline{\chi}\right)\Xi\left(t-\frac{iz}{2},\chi\right)y^{it}\, dt\nonumber\\
&\quad\quad\quad+\int_{0}^{\infty}f(z,t)\Xi\left(t+\frac{iz}{2},\overline{\chi}\right)\Xi\left(t-\frac{iz}{2},\chi\right)y^{-it}\, dt\bigg)\nonumber\\
&=\frac{1}{2}\bigg(\int_{0}^{\infty}f(z,t)\Xi\left(t+\frac{iz}{2},\overline{\chi}\right)\Xi\left(t-\frac{iz}{2},\chi\right)y^{it}\, dt\nonumber\\
&\quad\quad\quad+\int_{-\infty}^{0}f(z,-t)\Xi\left(-t+\frac{iz}{2},\overline{\chi}\right)\Xi\left(-t-\frac{iz}{2},\chi\right)y^{it}\, dt\bigg).\nonumber\\
\end{align}
However, using (\ref{fund}), we readily see that
\begin{align*}
\Xi\left(-t+\frac{iz}{2}, \overline{\chi}\right)&=\xi\left(\frac{1}{2}-it-\frac{z}{2},\overline{\chi}\right)=\epsilon(\chi)\xi\left(\frac{1}{2}+it+\frac{z}{2},\chi\right)=\epsilon(\chi)\Xi\left(t-\frac{iz}{2},\chi\right),\nonumber\\
\Xi\left(-t-\frac{iz}{2},\chi\right)&=\xi\left(\frac{1}{2}-it+\frac{z}{2},\chi\right)=\left(\epsilon(\chi)\right)^{-1}\xi\left(\frac{1}{2}+it-\frac{z}{2},\overline{\chi}\right)=\left(\epsilon(\chi)\right)^{-1}\Xi\left(t+\frac{iz}{2},\overline{\chi}\right),
\end{align*}
so that
\begin{equation}\label{ic1}
\Xi\left(-t+\frac{iz}{2}, \overline{\chi}\right)\Xi\left(-t-\frac{iz}{2},\chi\right)=\Xi\left(t+\frac{iz}{2},\overline{\chi}\right)\Xi\left(t-\frac{iz}{2},\chi\right).
\end{equation}
Thus from (\ref{ic}), (\ref{ic1}) and the fact that $f$ is an even function of $t$, we obtain
\begin{align*}
I(z,\mu,\chi)&=\frac{1}{2}\int_{-\infty}^{\infty}f(z,t)\Xi\left(t+\frac{iz}{2},\overline{\chi}\right)\Xi\left(t-\frac{iz}{2},\chi\right)y^{it}\, dt\\
&=\frac{1}{4i\sqrt{y}}\int_{\frac{1}{2}-i\infty}^{\frac{1}{2}+i\infty}\bigg(\phi\left(z,s-\frac{1}{2}\right)\phi\left(z,\frac{1}{2}-s\right)+\phi\left(-z,s-\frac{1}{2}\right)\phi\left(-z,\frac{1}{2}-s\right)\bigg)\nonumber\\
&\quad\quad\quad\quad\quad\quad\quad\times\xi\left(s-\frac{z}{2},\overline{\chi}\right)\xi\left(s+\frac{z}{2},\chi\right)y^{s}\, ds,
\end{align*}
where in the penultimate line, we made the change of variable $s=\frac{1}{2}+it$.
\end{proof}
For our purpose here, we replace $\mu$ by $2\mu$ in (\ref{genint2}) and then $t$ by $t/2$ on the left-hand side of (\ref{genint2}). Thus with $y=e^{2\mu}$, we find that
\begin{align}\label{gfica1}
&\int_{0}^{\infty}f\left(z,\frac{t}{2}\right)\Xi\left(\frac{t+iz}{2},\overline{\chi}\right)\Xi\left(\frac{t-iz}{2},\chi\right)\cos\mu t\, dt\nonumber\\
&=\frac{1}{2i\sqrt{y}}\int_{\frac{1}{2}-i\infty}^{\frac{1}{2}+i\infty}\bigg(\phi\left(z,s-\frac{1}{2}\right)\phi\left(z,\frac{1}{2}-s\right)+\phi\left(-z,s-\frac{1}{2}\right)\phi\left(-z,\frac{1}{2}-s\right)\bigg)\nonumber\\
&\quad\quad\quad\quad\quad\quad\quad\times\xi\left(s-\frac{z}{2},\overline{\chi}\right)\xi\left(s+\frac{z}{2},\chi\right)y^{s}\, ds.
\end{align}
It is this equation with which we will be working throughout this paper.
\section{Character analogues of the extended version of Guinand's formula}
We require the following lemma.
\begin{lemma}\label{lpro}
For \textup{Re} $s>1$ and \textup{Re} $(s-\eta)>1$, 
\begin{equation}\label{lprod}
L(s,\overline{\chi})L(s-\eta,\chi)=\sum_{n=1}^{\infty}\frac{\overline{\chi}(n)}{n^s}\sum_{d|n}\chi^{2}(d)d^\eta.
\end{equation}
\end{lemma}
\begin{proof}
Since the Dirichlet series for both the $L$-functions converge absolutely under the given hypotheses, using \cite[Theorem 11.5]{apostol}, we see that
\begin{align*}\label{prodl}
L(s,\overline{\chi})L(s-\eta,\chi)&=\sum_{n=1}^{\infty}\frac{\overline{\chi}(n)}{n^s}\sum_{k=1}^{\infty}\frac{\chi(k)}{k^{s-\eta}}\nonumber\\
&=\sum_{j=1}^{\infty}\frac{1}{j^s}\sum_{nk=j}\overline{\chi}(n)\chi(k)k^{\eta}\nonumber\\
&=\sum_{j=1}^{\infty}\frac{\overline{\chi}(j)}{j^s}\sum_{nk=j}\chi^{2}(k)k^{\eta},\nonumber\\
\end{align*}
since $\chi(k)\overline{\chi}(k)=1$.
\end{proof}
\begin{proof}[Theorem \textup{\ref{koshfec}}][]
First assume that $\chi$ is even. Let $\phi(z,s)\equiv 1$. Then from (\ref{finphi}), we see that $f(z,t)\equiv 1$. Using (\ref{Xic}), (\ref{xic}), (\ref{strivert}) and (\ref{lbound}), we find that the integral
\begin{equation*}
M(z, \mu, \chi):=\int_{0}^{\infty}\Xi\left(\frac{t+iz}{2},\overline{\chi}\right)\Xi\left(\frac{t-iz}{2},\chi\right)\cos\mu t\, dt
\end{equation*}
does converge. Using (\ref{gfica1}), we observe that
{\allowdisplaybreaks\begin{align}\label{gficake}
M(z, \mu, \chi)&=\frac{1}{i\sqrt{y}}\int_{\frac{1}{2}-i\infty}^{\frac{1}{2}+i\infty}\xi\left(s-\frac{z}{2},\overline{\chi}\right)\xi\left(s+\frac{z}{2},\chi\right)y^{s}\, ds\nonumber\\
&=\frac{1}{i\sqrt{y}}\int_{\frac{1}{2}-i\infty}^{\frac{1}{2}+i\infty}\Gamma\left(\frac{s}{2}-\frac{z}{4}\right)\Gamma\left(\frac{s}{2}+\frac{z}{4}\right)L\left(s-\frac{z}{2},\overline{\chi}\right)L\left(s+\frac{z}{2},\chi\right)\left(\frac{\pi}{qy}\right)^{-s}\, ds.
\end{align}}
Since Re $s=1/2$ and $-1<$ Re $z<1$, we have $0<$ Re $\left(s-\frac{z}{2}\right)<1$ and $0<$ Re $\left(s+\frac{z}{2}\right)<1$. Now replace $s$ by $s-\tf{z}{2}$ and let $\eta=-z$ in Lemma \ref{lpro}. Then, for Re $\left(s-\frac{z}{2}\right)>1$ and Re $\left(s+\frac{z}{2}\right)>1$,
\begin{equation}\label{usel}
L\left(s-\frac{z}{2},\overline{\chi}\right)L\left(s+\frac{z}{2},\chi\right)=\sum_{n=1}^{\infty}\frac{\overline{\chi}(n)}{n^{s-\frac{z}{2}}}\sum_{d|n}\chi^{2}(d)d^{-z}.
\end{equation}
We wish to shift the line of integration from Re $s=1/2$ to Re $s=3/2$ in order to be able to use (\ref{usel}) in (\ref{gficake}). Consider a positively oriented rectangular contour formed by $[\frac{1}{2}+iT, \frac{1}{2}-iT], [\frac{1}{2}-iT, \frac{3}{2}-iT], [\frac{3}{2}-iT,\frac{3}{2}+iT]$ and $[\frac{3}{2}+iT,\frac{1}{2}+iT]$, where $T$ is any positive real number. The integrand on the extreme right-hand side of (\ref{gficake}) does not have any pole inside the contour. Also as $T\to\infty$, the integrals along the horizontal segments $[\frac{1}{2}-iT, \frac{3}{2}-iT]$ and $[\frac{3}{2}+iT,\frac{1}{2}+iT]$ tend to zero, which can be seen by using (\ref{strivert}). Hence employing residue theorem, letting $T\to\infty$, using (\ref{usel}) in (\ref{gficake}), and interchanging the order of summation and integration because of absolute convergence, we observe that
\begin{align}\label{gficake1}
M(z, \mu, \chi)&=\frac{1}{i\sqrt{y}}\sum_{n=1}^{\infty}\overline{\chi}(n)n^{z/2}\bigg(\sum_{d|n}\chi^{2}(d)d^{-z}\bigg)\int_{\frac{3}{2}-i\infty}^{\frac{3}{2}+i\infty}\Gamma\left(\frac{s}{2}-\frac{z}{4}\right)\Gamma\left(\frac{s}{2}+\frac{z}{4}\right)\left(\frac{n\pi}{qy}\right)^{-s}\, ds.
\end{align}
But from \cite[p.~115, formula 11.1]{ober}, for $c=$ Re $s>\pm$ Re $\nu$,
\begin{equation}\label{bessmel}
\frac{1}{2\pi i}\int_{c-i\infty}^{c+i\infty}2^{s-2}w^{-s}\Gamma\left(\frac{s}{2}-\frac{\nu}{2}\right)\Gamma\left(\frac{s}{2}+\frac{\nu}{2}\right)x^{-s}\, ds=K_{\nu}(wx).
\end{equation}
Hence using (\ref{bessmel}) with $c=3/2$, $\nu=z/2$, $w=2$ and $x=n\pi/qy$ in (\ref{gficake1}), we find that
\begin{equation}\label{gficake2}
M(z, \mu, \chi)=\frac{8\pi}{\sqrt{y}}\sum_{n=1}^{\infty}\overline{\chi}(n)n^{z/2}\bigg(\sum_{d|n}\chi^{2}(d)d^{-z}\bigg)K_{\frac{z}{2}}\left(\frac{2\pi n}{qy}\right).
\end{equation}
Now let $\mu=\tf{1}{2}\log\alpha$ in (\ref{gficake2}) so that $y=e^{2\mu}$ implies that $y=\alpha$. Then using the fact that $\alpha\beta=1$, we deduce that
\begin{align}\label{gficake22}
&\frac{1}{8\pi}\int_{0}^{\infty}\Xi\left(\frac{t+iz}{2},\overline{\chi}\right)\Xi\left(\frac{t-iz}{2},\chi\right)\cos\left(\frac{1}{2} t\log\alpha\right)\, dt\nonumber\\
&=\sqrt{\beta}\sum_{n=1}^{\infty}\overline{\chi}(n)n^{z/2}\bigg(\sum_{d|n}\chi^{2}(d)d^{-z}\bigg)K_{\frac{z}{2}}\left(\frac{2\pi n\beta}{q}\right).
\end{align}
Next, observing that replacing $\alpha$ by $\beta$ and/or replacing simultaneously $\chi$ by $\overline{\chi}$ and $z$ by $-z$ in (\ref{gficake22}) leaves the integral on the left-hand side invariant, we obtain (\ref{koshfec1}).

Now consider the case when $\chi$ is odd. Again the convergence of the integral $M(z, \mu, \chi)$ can be seen from (\ref{strivert}) and (\ref{lbound}). Following similar steps above as in the case of even $\chi$, and using the definition of $\xi(s,\chi)$ from (\ref{xic}) for $\chi$ odd, we find that
\begin{equation}\label{gficako1}
M(z, \mu, \chi)=\frac{q}{i\pi\sqrt{y}}\sum_{n=1}^{\infty}\overline{\chi}(n)n^{z/2}\sum_{d|n}\chi^{2}(d)d^{-z}\int_{\frac{3}{2}-i\infty}^{\frac{3}{2}+i\infty}\Gamma\left(\frac{s}{2}-\frac{z}{4}+\frac{1}{2}\right)\Gamma\left(\frac{s}{2}+\frac{z}{4}+\frac{1}{2}\right)\left(\frac{n\pi}{qy}\right)^{-s}\, ds.
\end{equation}
Now replacing $s$ by $s+1$ in (\ref{bessmel}), we find that for $c=$ Re $s>\pm$ Re $\nu-1$,
\begin{equation}\label{bessmel1}
\frac{1}{2\pi i}\int_{c-i\infty}^{c+i\infty}2^{s-1}w^{-s-1}\Gamma\left(\frac{s+1}{2}-\frac{\nu}{2}\right)\Gamma\left(\frac{s+1}{2}+\frac{\nu}{2}\right)x^{-s}\, ds=xK_{\nu}(wx).
\end{equation}
Then using (\ref{bessmel1}) with $c=3/2$, $\nu=0$, $w=2$ and $x=n\pi/qy$ in (\ref{gficako1}), we see that
\begin{equation}\label{gficako2}
M(z, \mu, \chi)=\frac{8\pi}{y^{3/2}}\sum_{n=1}^{\infty}\overline{\chi}(n)n^{\frac{z}{2}+1}\bigg(\sum_{d|n}\chi^{2}(d)d^{-z}\bigg)K_{z/2}\left(\frac{2\pi n}{qy}\right).
\end{equation}
Now let $\mu=\tf{1}{2}\log\alpha$ in (\ref{gficako2}) so that $y=e^{2\mu}$ implies that $y=\alpha$. Then using the fact that $\alpha\beta=1$, we deduce that
\begin{align}\label{gficako22}
&\frac{1}{8\pi}\int_{0}^{\infty}\Xi\left(\frac{t+iz}{2},\overline{\chi}\right)\Xi\left(\frac{t-iz}{2},\chi\right)\cos\left(\frac{1}{2} t\log\alpha\right)\,dt\nonumber\\
&=\beta^{3/2}\sum_{n=1}^{\infty}\overline{\chi}(n)n^{\tf{z}{2}+1}\bigg(\sum_{d|n}\chi^{2}(d)d^{-z}\bigg)K_{z/2}\left(\frac{2\pi n\beta}{q}\right).
\end{align}
Next, observing that replacing $\alpha$ by $\beta$ and/or replacing simultaneously $\chi$ by $\overline{\chi}$ and $z$ by $-z$ in (\ref{gficako22}) leaves the integral on the left-hand side invariant, we obtain (\ref{koshfec1}).
\end{proof}
\textbf{Remark.} Letting $z\to 0$ in Theorem \ref{koshfec} gives a new character analogue of the extended version of Koshliakov's formula, i.e., Theorem \ref{koshfe2r}.\\

When $\chi$ is real, Theorem \ref{koshfec} reduces to the following corollary.
\begin{corollary}
Let $-1<$ Re $z<1$ and let $\chi$ denote a real, primitive, non-principal character modulo $q$. Let the number $b$ be defined as in \textup{(\ref{a})}. If
\begin{equation*}\label{fdefer}
F(z,\alpha,\chi)=\alpha^{b+\frac{1}{2}}\sum_{n=1}^{\infty}\chi(n)n^{-\frac{z}{2}+b}\sigma_{z}(n)K_{-\frac{z}{2}}\left(\frac{2\pi n\alpha}{q}\right),
\end{equation*}
then
\begin{align*}\label{koshfecr1}
F(z,\alpha,\chi)&=F(-z,\beta,\chi)=F(-z,\alpha,\chi)=F(z,\beta,\chi)\nonumber\\
&=\frac{1}{8\pi}\int_{0}^{\infty}\Xi\left(\frac{t+iz}{2},\chi\right)\Xi\left(\frac{t-iz}{2},\chi\right)\cos\left(\frac{1}{2}t\log\alpha\right)\, dt.
\end{align*}
\end{corollary}
The above corollary (without the integrals) is equivalent to the special cases, when $\chi$ is real, of the character analogues of Guinand's formula established in \cite{bds} (see Theorems 3.1 and 4.1). 
\section{Inverse Mellin transforms and asymptotic expansions of certain functions}
In this section, we evaluate inverse Mellin transforms of some functions and asymptotic expansions of certain other functions all of which are subsequently used in the later sections.
\begin{lemma}\label{invmelimp1}
For a primitive, non-principal character $\chi$, let $\psi(a,\chi)$ be defined as in \textup{(\ref{charpsi})}. Then for $0<c=$ \textup{Re} $s<1$ and $x\in\mathbb{R}\backslash\mathbb{Z}_{<0}$, 
\begin{equation}\label{invmelimp}
\frac{1}{2\pi i}\int_{c-i\infty}^{c+i\infty}\frac{L(1-s,\chi)}{\sin\pi s}x^{-s}\, ds=-\frac{1}{\pi}\psi(x,\chi).
\end{equation}
\end{lemma}
\begin{proof}
We will prove the result for even characters only. The case when $\chi$ is odd can be proved similarly. We first assume $|x|<1$ and later extend it to any real $x\in\mathbb{R}\backslash\mathbb{Z}_{<0}$ by analytic continuation. Let $0<c<1$. Consider a positively oriented rectangular contour formed by $[c-iT, c+iT], [c+iT, -M+iT], [-M+iT,-M-iT]$ and $[-M-iT,c-iT]$, where $T$ is any positive real number such that $T>\frac{\ln 2}{2\pi}$ and $M=n-1/2$ where $n$ is a positive integer. Let $s=\sigma+it$. Among the poles of the function $\left(L(1-s,\chi)/\sin\left(\pi s\right)\right) x^{-s}$, the only ones that contribute are the poles at the non-positive integers. Let $R_{a}$ denote the residue of the function $(L(1-s,\chi)/\sin\pi s) x^{-s}$ at $a$. Then,
\begin{equation}\label{rzero}
R_{0}=\lim_{s\to 0}\frac{sL(1-s,\chi)}{\sin\pi s}x^{-s}=\frac{1}{\pi}L(1,\chi).
\end{equation}
and
\begin{equation}\label{rn}
R_{-m}=\lim_{s\to -m}\frac{(s+m)L(1-s,\chi)}{\sin\pi s}x^{-s}=\frac{(-1)^{m}}{\pi}L(1+m,\chi)x^{m}.
\end{equation}
From (\ref{rzero}), (\ref{rn}) and the residue theorem, we have
\begin{align}\label{resap}
&\left[\int_{c-iT}^{c+iT}+\int_{c+iT}^{-M+iT}+\int_{-M+iT}^{-M-iT}+\int_{-M-iT}^{c-iT}\right]\frac{L(1-s,\chi)}{\sin\pi s}x^{-s}\, ds\nonumber\\
&=2\pi i\left(\frac{1}{\pi}L(1,\chi)+\sum_{0<m<M}\frac{(-1)^{m}}{\pi}L(1+m,\chi)x^{m}\right).
\end{align}
We first estimate the integrals along the upper and lower horizontal segments.
Using (\ref{lbound}), one finds that for $-M\leq\sigma\leq c$,
\begin{equation}\label{lbd}
L(1-\sigma\pm iT,\chi)=O\left(q^{c+1/2}T^{c+1/2}\right).
\end{equation}
Since $T>\frac{\ln 2}{2\pi}$, on the upper horizontal segment, we have 
\begin{equation}\label{uhssin}
\left|\frac{1}{\sin\pi s}\right|=\left|\frac{2e^{\pi is}}{e^{2\pi is}-1}\right|<4e^{-\pi T}.
\end{equation}
Similarly, on the lower horizontal segment,
\begin{equation*}\label{lhssin}
\left|\frac{1}{\sin\pi s}\right|=\left|\frac{2e^{-\pi is}}{1-e^{-2\pi is}}\right|<4e^{-\pi T}.
\end{equation*}
Since $|x|<1$, from (\ref{lbd}) and (\ref{uhssin}),
\begin{equation*}
\left|\int_{c+iT}^{-M+iT}\frac{L(1-s,\chi)}{\sin\pi s}x^{-s}\, ds\right|\leq
K_{1}(c+M)|x|^{-c}\cdot 4e^{-\pi T}q^{c+1/2}T^{c+1/2},
\end{equation*}
where $K_{1}$ is some constant.
Therefore,
\begin{equation}\label{uint}
\int_{c+i\infty}^{-M+i\infty}\frac{L(1-s,\chi)}{\sin\pi s}x^{-s}\, ds=0.
\end{equation}
Similarly,
\begin{equation}\label{lint}
\int_{-M-i\infty}^{c-i\infty}\frac{L(1-s,\chi)}{\sin\pi s}x^{-s}\, ds=0.
\end{equation}
Then from (\ref{resap}), (\ref{uint}) and (\ref{lint}), we find that
\begin{equation}\label{resap1}
\left[\int_{c-i\infty}^{c+i\infty}+\int_{-M+i\infty}^{-M-i\infty}\right]\frac{L(1-s,\chi)}{\sin\pi s}x^{-s}\, ds
=2\pi i\left(\frac{1}{\pi}L(1,\chi)+\sum_{0<m<M}\frac{(-1)^{m}}{\pi}L(1+m,\chi)x^{m}\right).
\end{equation}
It remains to examine $\int_{-M+i\infty}^{-M-i\infty}(L(1-s,\chi)/\sin\pi s)x^{-s}\, ds$. Since $M=n-1/2$, we have
$|\sin\pi(-M+it)|=|\cosh\pi t|\geq 1$ and $L(1+M\pm it)=O(1)$ as $1+M>1$. Thus,
\begin{align*}
\left|\int_{-M-i\infty}^{-M+i\infty}\frac{L(1-s,\chi)}{\sin\pi s}x^{-s}\, ds\right|&=\left|i\int_{-\infty}^{\infty}\frac{L(1+M-it,\chi)}{\sin\pi (-M+it)}x^{M-it}\, dt\right|\nonumber\\
&=|x|^{M}\int_{-1}^{1}O(1)\, dt+|x|^{M}\int_{1}^{\infty}O\left(e^{-\pi|t|}\right)\, dt+|x|^{M}\int_{-\infty}^{-1}O\left(e^{-\pi|t|}\right)\, dt\nonumber\\
&=O(|x|^{M}),
\end{align*}
as $|t|\to\infty$.
Since $|x|<1$,
\begin{equation}\label{vint}
\lim_{M\to\infty}\int_{-M+i\infty}^{-M-i\infty}\frac{L(1-s,\chi)}{\sin\pi s}x^{-s}\, ds=0.
\end{equation}
From (\ref{resap1}) and (\ref{vint}), we see that
\begin{align}
\frac{1}{2\pi i}\int_{c-i\infty}^{c+i\infty}\frac{L(1-s,\chi)}{\sin\pi s}x^{-s}\,ds &=\left(\frac{1}{\pi}L(1,\chi)+\sum_{m=1}^{\infty}\frac{(-1)^{m}}{\pi}L(1+m,\chi)x^{m}\right)\nonumber\\
&=\frac{-1}{\pi}\left(-L(1,\chi)-\sum_{k=1}^{\infty}\frac{\chi(k)}{k}\sum_{m=1}^{\infty}\left(\frac{-x}{k}\right)^{m}\right)\nonumber\\
&=\frac{-1}{\pi}\left(-L(1,\chi)-\sum_{k=1}^{\infty}\chi(k)\left(\frac{1}{x+k}-\frac{1}{k}\right)\right)\nonumber\\
&=\frac{-1}{\pi}\psi(x,\chi).
\end{align}
Since both sides of (\ref{invmelimp}) are analytic for any $x\in\mathbb{R}\backslash\mathbb{Z}_{<0}$, the result follows by analytic continuation.
\end{proof}
\begin{lemma}\label{invmelimpeg}
Let $z\in\mathbb{C}$ be fixed such that $-1<$ \textup{Re} $z<1$. For a primitive, non-principal character $\chi$, let $L(z, a, \chi)$ be defined as in \textup{(\ref{charhur})}. Then for $-\frac{1}{2}$ \textup{Re} $z<c=$ \textup{Re} $s<\frac{1}{2}$ \textup{Re} $z$ and $x\in\mathbb{R}\backslash\mathbb{Z}_{<0}$, 
\begin{equation}\label{invmelimpg}
\frac{1}{2\pi i}\int_{c-i\infty}^{c+i\infty}\Gamma\left(s+\frac{z}{2}\right)\Gamma\left(1-s+\frac{z}{2}\right)L\left(1-s+\frac{z}{2},\chi\right)x^{-s}\, ds=x^{z/2}\Gamma(z+1)L(z+1,x,\chi).
\end{equation}
\end{lemma}
\begin{proof}
We prove the result only for even characters. The case for odd characters can be proved similarly. We first assume $|x|<1$ and later extend it to any real $x\in\mathbb{R}\backslash\mathbb{Z}_{<0}$ by analytic continuation. Let $-\frac{1}{2}$ \textup{Re} $z<c=$ \textup{Re} $s<\frac{1}{2}$ \textup{Re} $z$. Consider a positively oriented rectangular contour formed by $[c-iT, c+iT], [c+iT, -M+iT], [-M+iT,-M-iT]$ and $[-M-iT,c-iT]$, where $T$ is some positive real number and $M=n-1/2$, where $n$ is a positive integer. Let $s=\sigma+it$. Among the poles of the function $\displaystyle\Gamma(s+z/2)\Gamma(1-s+z/2)L(1-s+z/2,\chi)x^{-s}$, the only ones that contribute are the poles at $s=-z/2-m, m\geq 0$. Let $R_{a}$ denote the residue of the function $\displaystyle\Gamma(s+z/2)\Gamma(1-s+z/2)L(1-s+z/2,\chi)x^{-s}$ at $a$. Then for $m\geq 0$,
\begin{align}\label{resneww}
R_{-\frac{z}{2}-m}&=\lim_{s\to-\frac{z}{2}-m}\left(s+\frac{z}{2}+m\right)\Gamma\left(s+\frac{z}{2}\right)\Gamma\left(1-s+\frac{z}{2}\right)L\left(1-s+\frac{z}{2},\chi\right)x^{-s}\nonumber\\
&=\frac{(-1)^{m}}{m!}\Gamma\left(1+z+m\right)L\left(1+z+m,\chi\right)x^{z/2+m}.
\end{align}
From (\ref{resneww}) and the residue theorem, we have
\begin{align}\label{resapw}
&\left[\int_{c-iT}^{c+iT}+\int_{c+iT}^{-M+iT}+\int_{-M+iT}^{-M-iT}+\int_{-M-iT}^{c-iT}\right]\Gamma\left(s+\frac{z}{2}\right)\Gamma\left(1-s+\frac{z}{2}\right)L\left(1-s+\frac{z}{2},\chi\right)x^{-s}\, ds\nonumber\\
&=2\pi ix^{z/2}\sum_{0\leq m<M}\frac{(-1)^{m}}{m!}\Gamma\left(1+z+m\right)L\left(1+z+m,\chi\right)x^{m}.
\end{align}
We now estimate the integral along the upper horizontal segment. From (\ref{lbd}), we easily see that for $-M\leq\sigma\leq c$, i.e., $-M-$ Re $\frac{z}{2}\leq\sigma-$ Re $\frac{z}{2}\leq c-$ Re $\frac{z}{2}$,
\begin{equation}\label{lbdw}
L\left(1-\left(\sigma-\textup{Re}\frac{z}{2}\right)-i\left(T-\textup{Im}\frac{z}{2}\right),\chi\right)=O\left(q^{c-\textup{Re}\frac{z}{2}+\frac{1}{2}}\left(T-\textup{Im}\frac{z}{2}\right)^{c-\textup{Re}\frac{z}{2}+\frac{1}{2}}\right).
\end{equation}
By (\ref{strivert}), we observe that
\begin{equation}\label{g1}
\left|\Gamma\left(s+\frac{z}{2}\right)\right|\sim\sqrt{2\pi}e^{-\frac{\pi}{2}\left|T+\textup{Im}\frac{z}{2}\right|}\cdot\left|T+\textup{Im}\frac{z}{2}\right|^{\sigma+\textup{Re}\frac{z}{2}-\frac{1}{2}},
\end{equation}
and
\begin{equation}\label{g2}
\left|\Gamma\left(1-s+\frac{z}{2}\right)\right|\sim\sqrt{2\pi}e^{-\frac{\pi}{2}\left|T-\textup{Im}\frac{z}{2}\right|}\cdot\left|T-\textup{Im}\frac{z}{2}\right|^{-\sigma+\textup{Re}\frac{z}{2}+\frac{1}{2}},
\end{equation}
Since $|x|<1$, from (\ref{lbdw}), (\ref{g1}) and (\ref{g2}), we deduce that
\begin{align*}\label{uhsw}
&\left|\int_{c+iT}^{-M+iT}\Gamma(s+z/2)\Gamma(1-s+z/2)L(1-s+z/2,\chi)x^{-s}\, ds\right|\nonumber\\
&\leq 2\pi K_{3}(c+M)|x|^{-c}e^{-\frac{\pi}{2}\left(\left|T+\textup{Im}\frac{z}{2}\right|+\left|T-\textup{Im}\frac{z}{2}\right|\right)}\left|T+\textup{Im}\frac{z}{2}\right|^{\sigma+\textup{Re}\frac{z}{2}-\frac{1}{2}}\left|T-\textup{Im}\frac{z}{2}\right|^{-\sigma+\textup{Re}\frac{z}{2}+\frac{1}{2}},
\end{align*}
where $K_{3}$ is a constant. Hence, 
\begin{equation}\label{uhswf}
\int_{c+i\infty}^{-M+i\infty}\Gamma(s+z/2)\Gamma(1-s+z/2)L(1-s+z/2,\chi)x^{-s}\, ds=0.
\end{equation}
Similarly for the integral along the lower horizontal segment, using (\ref{g1}), (\ref{g2}) and the fact that 
\begin{equation*}\label{lbdw2}
L\left(1-\left(\sigma-\textup{Re}\frac{z}{2}\right)+i\left(T+\textup{Im}\frac{z}{2}\right),\chi\right)=O\left(q^{c-\textup{Re}\frac{z}{2}+\frac{1}{2}}\left(T+\textup{Im}\frac{z}{2}\right)^{c-\textup{Re}\frac{z}{2}+\frac{1}{2}}\right),
\end{equation*}
we observe that
\begin{equation}\label{lhswf}
\int_{-M-i\infty}^{c-i\infty}\Gamma(s+z/2)\Gamma(1-s+z/2)L(1-s+z/2,\chi)x^{-s}\, ds=0.
\end{equation}
Hence from (\ref{resapw}), (\ref{uhswf}) and (\ref{lhswf}), it is clear that
\begin{align}\label{sti2}
&\left[\int_{c-i\infty}^{c+i\infty}+\int_{-M+i\infty}^{-M-i\infty}\right]\Gamma(s+z/2)\Gamma(1-s+z/2)L(1-s+z/2,\chi)x^{-s}\, ds\nonumber\\
&=2\pi ix^{z/2}\sum_{0\leq m<M}\frac{(-1)^{m}}{m!}\Gamma\left(1+z+m\right)L\left(1+z+m,\chi\right)x^{m}.
\end{align}
It remains to evaluate $\int_{-M+i\infty}^{-M-i\infty}\Gamma(s+z/2)\Gamma(1-s+z/2)L(1-s+z/2,\chi)x^{-s}\, ds$. In \cite{bs}, we find that as $|t|\to\infty$,
\begin{equation*}
\Gamma\left(-M+it\right)=O\left(|t|^{-M-1/2}e^{-\pi|t|/2}\right).
\end{equation*}
Hence as $T\to\infty$, 
\begin{equation}\label{fg1}
\Gamma\left(-M+it+\frac{z}{2}\right)=O\left(\left|T+\textup{Im}\frac{z}{2}\right|^{-M+\textup{Re}\frac{z}{2}-\frac{1}{2}}e^{-\frac{\pi}{2}|T+\textup{Im}\frac{z}{2}|}\right).
\end{equation}
Again by (\ref{strivert}), as $T\to\infty$,
\begin{equation}
\left|\Gamma\left(1+M-it+\frac{z}{2}\right)\right|=\sqrt{2\pi}e^{-\frac{\pi}{2}\left|T-\textup{Im}\frac{z}{2}\right|}\cdot\left|T-\textup{Im}\frac{z}{2}\right|^{M+\textup{Re}\frac{z}{2}+\frac{1}{2}}\left(1+O\left(\frac{1}{\left|T-\textup{Im}\frac{z}{2}\right|}\right)\right).
\end{equation}
Also, $L\left(1+M-it+\frac{z}{2},\chi\right)$ is bounded as Re$\left(1+M-it+\frac{z}{2}\right)>1$. Hence,
\begin{align*}
&\left|\int_{-M+i\infty}^{-M-i\infty}\Gamma\left(s+\frac{z}{2}\right)\Gamma\left(1-s+\frac{z}{2}\right)L\left(1-s+\frac{z}{2},\chi\right)x^{-s}\, ds\right|\nonumber\\
&=\left|i\int_{-\infty}^{\infty}\Gamma\left(-M+it+\frac{z}{2}\right)\Gamma\left(1+M-it+\frac{z}{2}\right)L\left(1+M-it+\frac{z}{2},\chi\right)x^{M-it}\, dt\right|\nonumber\\
&=|x|^{M}\int_{-1}^{1}O(1)\, dt+|x|^{M}\int_{1}^{\pm\infty}O\left(\left|T+\textup{Im}\frac{z}{2}\right|^{-M+\textup{Re}\frac{z}{2}-\frac{1}{2}}\left|T-\textup{Im}\frac{z}{2}\right|^{M+\textup{Re}\frac{z}{2}+\frac{1}{2}}e^{-\frac{\pi}{2}\left(\left|T+\textup{Im}\frac{z}{2}\right|+\left|T-\textup{Im}\frac{z}{2}\right|\right)}\right)dt\nonumber\\
&=O\left(|x|^{M}\right),
\end{align*}
as $T\to\infty$. Since $|x|<1$,
\begin{equation}\label{wow}
\lim_{M\to\infty}\int_{-M+i\infty}^{-M-i\infty}\Gamma\left(s+\frac{z}{2}\right)\Gamma\left(1-s+\frac{z}{2}\right)L\left(1-s+\frac{z}{2},\chi\right)x^{-s}\, ds=0.
\end{equation}
From (\ref{sti2}), (\ref{wow}), we finally deduce that
{\allowdisplaybreaks\begin{align*}
&\frac{1}{2\pi i}\int_{c-i\infty}^{c+i\infty}\Gamma\left(s+\frac{z}{2}\right)\Gamma\left(1-s+\frac{z}{2}\right)L\left(1-s+\frac{z}{2},\chi\right)x^{-s}\, ds\nonumber\\
&=x^{z/2}\sum_{m=0}^{\infty}\frac{(-1)^{m}}{m!}\Gamma\left(1+z+m\right)L\left(1+z+m,\chi\right)x^{z/2+m}\nonumber\\
&=x^{z/2}\Gamma(z+1)\sum_{m=0}^{\infty}\frac{(-1)^{m}}{m!}\frac{\Gamma\left(1+z+m\right)}{\Gamma(1+z)}\sum_{k=1}^{\infty}\frac{\chi(k)}{k^{z+m+1}}x^{m}\nonumber\\
&=x^{z/2}\Gamma(z+1)\sum_{k=1}^{\infty}\frac{\chi(k)}{k^{z+1}}\sum_{m=0}^{\infty}\frac{\Gamma(1+z+m)}{m!\Gamma(1+z)}\left(\frac{-x}{k}\right)^{m}\nonumber\\
&=x^{z/2}\Gamma(z+1)\sum_{k=1}^{\infty}\frac{\chi(k)}{k^{z+1}}\left(1+\frac{x}{k}\right)^{-z-1}\nonumber\\
&=x^{z/2}\Gamma(z+1)\sum_{k=1}^{\infty}\frac{\chi(k)}{(k+x)^{z+1}}\nonumber\\
&=x^{z/2}\Gamma(z+1)L(z+1,x,\chi),
\end{align*}}%
where in the fourth step above, we have utilized binomial theorem since $|x|<1$. Since both sides of (\ref{invmelimpg}) are analytic for any $x\in\mathbb{R}\backslash\mathbb{Z}_{<0}$, the result follows by analytic continuation.
\end{proof}
For $j\geq 1$, the generalized Bernoulli numbers $B_{j}(\chi)$ are given by \cite[p.~426]{berndt}
\begin{equation*}\label{bere}
B_{2j}(\chi)=\frac{2(-1)^{j-1}G(\overline{\chi})(2j)!}{k(2\pi/k)^{2j}}L(2j,\chi),
\end{equation*}
for $\chi$ even, and by
\begin{equation*}\label{bero}
B_{2j-1}(\chi)=\frac{2(-1)^{j-1}iG(\overline{\chi})(2j-1)!}{k(2\pi/k)^{2j-1}}L(2j-1,\chi),
\end{equation*}
for $\chi$ odd. Also it is known \cite[p.~423, Corollary 3.4]{berndt} that $B_{2j-1}(\chi)=0$ when $\chi$ is even and $B_{2j}(\chi)=0$ when $\chi$ is odd. Next, we give an asymptotic expansion of $\psi(a,\chi)$ as $|a|\to\infty$.
\begin{lemma}\label{asycharpsi}
For $-\pi<\arg a<\pi$, as $|a|\to\infty$,
\begin{equation*}\label{acp}
\psi(a,\chi)\sim-\frac{L(0,\chi)}{a}-\chi(-1)\sum_{j=2}^{\infty}\frac{B_{j}(\overline{\chi})}{ja^{j}}.
\end{equation*}
\end{lemma}
\begin{proof}
The proof is very much identical to Proposition 5.3 in \cite[p.~435]{berndt} and hence is omitted.
\end{proof}
The asymptotic expansion of $L(z,a,\chi)$ as $|a|\to\infty$ is given below.
\begin{lemma}\label{asympg}
For \textup{Re} $z>0$ and $-\pi<\arg a<\pi$, as $|a|\to\infty$,
\begin{equation*}\label{acpg}
L(z,a,\chi)\sim\chi(-1)\sum_{j=1}^{\infty}\frac{B_{j}(\overline{\chi})\prod_{m=0}^{j-2}(z+m)}{j!a^{z+j-1}}.
\end{equation*}
\end{lemma}
\begin{proof}
To derive this, one takes (4.3) and (4.4) in \cite[p.~424]{berndt} valid for $\chi$ even and odd respectively, substitutes $A=0, B=N, r=1$ and $f(u)=(u+a)^{-z}$, lets $N\to\infty$ and then performs repeated integration by parts on the prevalent integral.
\end{proof}
\section{Character analogues of Theorem \ref{sechur}}
In this section, we prove analogues of Theorem \ref{sechur} for even and odd primitive characters. Then we give character analogues of Ramanujan's transformation formula (Theorem \ref{entry1}) as special cases.
\begin{proof}[Theorem \textup{\ref{ramgenec}}][]
Using Lemma \ref{asympg}, one sees that the series involving the functions $L(z,a,\chi)$ in the theorem are convergent. Let $\phi(z,s)=(z+1+2s)\Gamma\left(\frac{-z-1}{4}+\frac{s}{2}\right)$. Then from (\ref{finphi}) and (\ref{omdef}), we find that $f\left(z,\frac{t}{2}\right)=\frac{1}{2}\Omega(z,t)$. From (\ref{gfica1}), we have
\begin{align}\label{rockon}
\int_{0}^{\infty}\Omega(z,t)\Xi\left(\frac{t+iz}{2},\overline{\chi}\right)\Xi\left(\frac{t-iz}{2},\chi\right)\cos\mu t\, dt=\frac{1}{i\sqrt{y}}\left(J(z,y,\chi)+J(-z,y,\overline{\chi})\right),
\end{align}
where
\begin{equation}\label{intt}
J(z,y,\chi):=\int_{\frac{1}{2}-i\infty}^{\frac{1}{2}+i\infty}U(z,s,y,\chi)\, ds
\end{equation}
with
\begin{equation*}\label{intu}
U(z,s,y,\chi):=(-z+2s)(-z+2-2s)\Gamma\left(\frac{z}{4}+\frac{s}{2}-\frac{1}{2}\right)\Gamma\left(\frac{z}{4}-\frac{s}{2}\right)\xi\left(s-\frac{z}{2},\overline{\chi}\right)\xi\left(s+\frac{z}{2},\chi\right)y^{s}.
\end{equation*}
Using (\ref{strivert}) and (\ref{lbound}), one sees that indeed the integral on the left side of (\ref{rockon}) converges. We first simplify the integrand in (\ref{intt}). Using (\ref{xic}) with $b=0$, and then duplication formula \cite[p.~46, Equation (3.4)]{temme} and reflection formula \cite[p.~46, Equation (3.5)]{temme} for Gamma function in the second equality below, we have
\begin{align}\label{u}
U(z,s,y,\chi)&=16\left(\frac{\pi}{qy}\right)^{-s}\left\{\Gamma\left(\frac{z}{4}+\frac{s+1}{2}\right)\Gamma\left(\frac{z}{4}+\frac{s}{2}\right)\right\}\left\{\Gamma\left(\frac{z}{4}-\frac{s}{2}+1\right)\Gamma\left(\frac{s}{2}-\frac{z}{4}\right)\right\}\nonumber\\
&\quad\times L\left(s-\frac{z}{2},\overline{\chi}\right)L\left(s+\frac{z}{2},\chi\right)\nonumber\\
&=16\left(\frac{\pi}{qy}\right)^{-s}\cdot\frac{\sqrt{\pi}}{2^{s+\frac{z}{2}-1}}\Gamma\left(s+\frac{z}{2}\right)\cdot \frac{\pi}{\sin\left(\pi\left(\frac{s}{2}-\frac{z}{4}\right)\right)}\cdot L\left(s-\frac{z}{2},\overline{\chi}\right)L\left(s+\frac{z}{2},\chi\right).
\end{align}
Substituting (\ref{funl}) in the form
\begin{equation*}\label{funld}
L\left(s-\frac{z}{2},\overline{\chi}\right)=\frac{(2\pi)^{s-\frac{z}{2}}L(1-s+\frac{z}{2},\chi)}{2q^{s-\frac{z}{2}-1}G(\chi)\Gamma\left(s-\frac{z}{2}\right)\cos\left(\frac{\pi}{2}\left(s-\frac{z}{2}\right)\right)}
\end{equation*}
in (\ref{u}) and then simplifying, we find that
\begin{equation}\label{u1}
U(z,s,y,\chi)=\frac{32y^{s}2^{-z}\pi^{(1-z)/2}}{q^{-\frac{z}{2}-1}G(\chi)}\Gamma\left(1-s+\frac{z}{2}\right)\Gamma\left(s+\frac{z}{2}\right)L\left(1-s+\frac{z}{2},\chi\right)L\left(s+\frac{z}{2},\chi\right).
\end{equation}
We wish to shift the line of integration from Re $s=1/2$ to Re $s=3/2$ in order to evaluate (\ref{intt}), since then $-1<$ Re $z<1$ implies that Re $(s+z/2)>1$, which allows us to use the series representation of $L\left(s+\frac{z}{2},\chi\right)$. Consider a positively oriented rectangular contour formed by $[\frac{1}{2}+iT, \frac{1}{2}-iT], [\frac{1}{2}-iT, \frac{3}{2}-iT], [\frac{3}{2}-iT,\frac{3}{2}+iT]$ and $[\frac{3}{2}+iT,\frac{1}{2}+iT]$, where $T$ is any positive real number. The integrand in (\ref{intt}) does not have any pole inside the contour since the pole of $\Gamma\left(1-s+\frac{z}{2}\right)$ at $s=1+z/2$ is cancelled by the zero of $L\left(1-s+\frac{z}{2},\chi\right)$ there. Also as $T\to\infty$, the integrals along the horizontal segments $[\frac{1}{2}-iT, \frac{3}{2}-iT]$ and $[\frac{3}{2}+iT,\frac{1}{2}+iT]$ tend to zero, which can be seen using (\ref{strivert}). Employing residue theorem, letting $T\to\infty$ and using (\ref{u1}), we find that
\begin{align}\label{intt1}
J(z,y,\chi)
&=\frac{32\cdot 2^{-z}\pi^{(1-z)/2}}{q^{-\frac{z}{2}-1}G(\chi)}\int_{\frac{3}{2}-i\infty}^{\frac{3}{2}+i\infty}\Gamma\left(s+\frac{z}{2}\right)\Gamma\left(1-s+\frac{z}{2}\right)L\left(1-s+\frac{z}{2},\chi\right)L\left(s+\frac{z}{2},\chi\right)y^{s}\, ds\nonumber\\
&=\frac{32\cdot 2^{-z}\pi^{(1-z)/2}}{q^{-\frac{z}{2}-1}G(\chi)}\sum_{n=1}^{\infty}\frac{\chi(n)}{n^{z/2}}\int_{\frac{3}{2}-i\infty}^{\frac{3}{2}+i\infty}\Gamma\left(s+\frac{z}{2}\right)\Gamma\left(1-s+\frac{z}{2}\right)L\left(1-s+\frac{z}{2},\chi\right)\left(\frac{n}{y}\right)^{-s}\, ds.
\end{align}
Now, in order to use Lemma \ref{invmelimpeg} to evaluate the integral in (\ref{intt1}), we again have to shift the line of integration from Re $s>3/2$ to Re $s=d$, where $-\frac{1}{2}$ \textup{Re} $z<d<\frac{1}{2}$ \textup{Re} $z$. Again, we do not encounter any pole in this process. Hence,
\begin{equation}\label{intt3}
J(z,y,\chi)=\frac{64i2^{-z}y^{-z/2}\pi^{(3-z)/2}\Gamma(z+1)}{q^{-\frac{z}{2}-1}G(\chi)}\sum_{n=1}^{\infty}\sum_{k=1}^{\infty}\frac{\chi(n)\chi(k)}{(k+n/y)^{z+1}}.
\end{equation}
Since $-1<\textup{Re}(z)<1$, the other integral, namely $J(-z,y,\overline{\chi})$, can be evaluated by simply replacing $z$ by $-z$ and $\chi$ by $\overline{\chi}$ in (\ref{intt3}). Now (\ref{rockon}), (\ref{intt3}), (\ref{te}) and the discussion in the previous line give
{\allowdisplaybreaks\begin{align}\label{pt}
&\int_{0}^{\infty}\Omega(z,t)\Xi\left(\frac{t+iz}{2},\overline{\chi}\right)\Xi\left(\frac{t-iz}{2},\chi\right)\cos\mu t\, dt\nonumber\\
&=\frac{64\pi^{3/2}q}{\sqrt{y}}\left(T(z,y^{-1},\chi)\sum_{n=1}^{\infty}\sum_{k=1}^{\infty}\frac{\chi(n)\chi(k)}{(k+n/y)^{z+1}}+T(-z,y^{-1},\overline{\chi})\sum_{n=1}^{\infty}\sum_{k=1}^{\infty}\frac{\overline{\chi}(n)\overline{\chi}(k)}{(k+n/y)^{-z+1}}\right),
\end{align}}
where it is easy to see from the fact that $-1<$ Re $z<1$, from the discussion just preceding the statement of Theorem \ref{ramgenec} and from Lemma \ref{asympg} that both the double series on the right-hand side of (\ref{pt}) converge.

Now let $\mu=\tf{1}{2}\log\alpha$ in (\ref{pt}) so that $y=e^{2\mu}$ implies that $y=\alpha$. Then using the fact that $\alpha\beta=1$ and using (\ref{charhur}) in the second equality below, we deduce that
{\allowdisplaybreaks\begin{align*}\label{pt1}
&\int_{0}^{\infty}\Omega(z,t)\Xi\left(\frac{t+iz}{2},\overline{\chi}\right)\Xi\left(\frac{t-iz}{2},\chi\right)\cos\left(\tf{1}{2}t\log\alpha\right)\, dt\nonumber\\
&=64\pi^{3/2}q\sqrt{\beta}\left(T(z,\beta,\chi)\sum_{n=1}^{\infty}\sum_{k=1}^{\infty}\frac{\chi(n)\chi(k)}{(k+n\beta)^{z+1}}+T(-z,\beta,\overline{\chi})\sum_{n=1}^{\infty}\sum_{k=1}^{\infty}\frac{\overline{\chi}(n)\overline{\chi}(k)}{(k+n\beta)^{-z+1}}\right)\nonumber\\
&=64\pi^{3/2}q\sqrt{\beta}\left(T(z,\beta,\chi)\sum_{n=1}^{\infty}\chi(n)L(z+1,n\beta,\chi)+T(-z,\beta,\overline{\chi})\sum_{n=1}^{\infty}\overline{\chi}(n)L(-z+1,n\beta,\overline{\chi})\right).
\end{align*}}
The integral on the extreme left-hand side above is invariant under the transformation $\alpha\to\beta$ or under the simultaneous application of the transformations $\alpha\to\beta, \chi\to\overline{\chi}$ and $z\to -z$. Thus we obtain (\ref{ramgenec1}).
\end{proof}
Next we give an analogue of Ramanujan's transformation formula (Theorem \ref{entry1}) for even characters.
\begin{corollary}\label{1splcases} 
For an even character $\chi$, define $P(\alpha,\chi)$ by
\begin{equation*}\label{pdefi}
P(\alpha,\chi) = \sqrt{\alpha}\hspace{1mm} \textup{Re} \left(G(\chi)\sum_{n=1}^{\infty}\sum_{k=1}^{\infty}\frac{\overline{\chi}(n)\overline{\chi}(k)}{k+n\alpha}\right) =-\sqrt{\alpha}\hspace{1mm} \textup{Re}\left(G(\chi)\sum_{n=1}^{\infty}\overline{\chi}(n)\psi\left(n\alpha,\overline{\chi}\right)\right),
\end{equation*}
where $\psi\left(a,\chi\right)$ is defined in \textup{(\ref{charpsi})}. Then we have
\begin{align}\label{splcases}
&P(\alpha,\chi)=P(\beta,\overline{\chi})=P(\alpha,\overline{\chi})=P(\beta,\chi)\nonumber\\
&=\frac{1}{64\pi^{3/2}}\int_{0}^{\infty}(1+t^2)\Gamma\bigg(\frac{-1+it}{4}\bigg)\Gamma\bigg(\frac{-1-it}{4}\bigg)\Xi\left(\frac{t}{2},\overline{\chi}\right)\Xi\left(\frac{t}{2},\chi\right)\cos\left(\frac{1}{2}t\log\alpha\right)\, dt.
\end{align}
\end{corollary}
\begin{proof}
Using Lemma \ref{asycharpsi}, we readily see that the double series in the definition of $P(\alpha,\chi)$ converges. Let $z\to 0$ in (\ref{ramgenec1}). Then multiplying both sides by $q$ and using (\ref{eo}), we have
{\allowdisplaybreaks\begin{align}\label{splcase}
&\sqrt{\alpha}\left(G(\overline{\chi})\sum_{n=1}^{\infty}\sum_{k=1}^{\infty}\frac{\chi(n)\chi(k)}{k+n\alpha}+G(\chi)\sum_{n=1}^{\infty}\sum_{k=1}^{\infty}\frac{\overline{\chi}(n)\overline{\chi}(k)}{k+n\alpha}\right)\nonumber\\
&=\sqrt{\beta}\left(G(\chi)\sum_{n=1}^{\infty}\sum_{k=1}^{\infty}\frac{\overline{\chi}(n)\overline{\chi}(k)}{k+n\beta}+G(\overline{\chi})\sum_{n=1}^{\infty}\sum_{k=1}^{\infty}\frac{\chi(n)\chi(k)}{k+n\beta}\right)\nonumber\\
&=\frac{1}{32\pi^{3/2}}\int_{0}^{\infty}(1+t^2)\Gamma\bigg(\frac{-1+it}{4}\bigg)\Gamma\bigg(\frac{-1-it}{4}\bigg)\Xi\left(\frac{t}{2},\overline{\chi}\right)\Xi\left(\frac{t}{2},\chi\right)\cos\left(\frac{1}{2}t\log\alpha\right)\, dt.
\end{align}}
Each of the first two expressions in (\ref{splcase}) can be written in two different ways as real parts of a double series. Thus,
{\allowdisplaybreaks\begin{align*}\label{innov}
&\sqrt{\alpha}\hspace{1mm}\textup{Re} \left(G(\overline{\chi})\sum_{n=1}^{\infty}\sum_{k=1}^{\infty}\frac{\chi(n)\chi(k)}{k+n\alpha}\right)=\sqrt{\alpha}\hspace{1mm} \textup{Re} \left(G(\chi)\sum_{n=1}^{\infty}\sum_{k=1}^{\infty}\frac{\overline{\chi}(n)\overline{\chi}(k)}{k+n\alpha}\right)\nonumber\\
&=\sqrt{\beta}\hspace{1mm}\textup{Re}\left(G(\chi)\sum_{n=1}^{\infty}\sum_{k=1}^{\infty}\frac{\overline{\chi}(n)\overline{\chi}(k)}{k+n\beta}\right)=\sqrt{\beta}\hspace{1mm} \textup{Re} \left(G(\overline{\chi})\sum_{n=1}^{\infty}\sum_{k=1}^{\infty}\frac{\chi(n)\chi(k)}{k+n\beta}\right)\nonumber\\
&=\frac{1}{64\pi^{3/2}}\int_{0}^{\infty}(1+t^2)\Gamma\bigg(\frac{-1+it}{4}\bigg)\Gamma\bigg(\frac{-1-it}{4}\bigg)\Xi\left(\frac{t}{2},\overline{\chi}\right)\Xi\left(\frac{t}{2},\chi\right)\cos\left(\frac{1}{2}t\log\alpha\right)\, dt.
\end{align*}}
This implies (\ref{splcases}). 
\end{proof}
Moreover, if we start with the integral in Corollary \ref{1splcases}, evaluate it using (\ref{gfica1}) with $z=0$ and make use of Lemma \ref{invmelimp1} when $\chi$ is even, we obtain the same result as in Corollary \ref{1splcases}, except that the function $P(\alpha,\chi)$ is replaced by the function $F(\alpha,\chi)$ defined by
\begin{equation}\label{ta}
F(\alpha,\chi):=\sqrt{\alpha}G(\chi)\sum_{n=1}^{\infty}\sum_{k=1}^{\infty}\frac{\overline{\chi}(n)\overline{\chi}(k)}{k+n\alpha}=-\sqrt{\alpha}G(\chi)\sum_{n=1}^{\infty}\overline{\chi}(n)\psi\left(n\alpha,\overline{\chi}\right).
\end{equation}
It is then trivial to see that $F(\alpha,\chi)=P(\alpha,\chi)$. 

Theorem \ref{ramgenoc} can be analogously proved using Lemma \ref{invmelimpeg} for $\chi$ odd. We just note that there we have to take care of the pole of $\Gamma\left(1-s+\frac{1}{2}z\right)$ in the integrands of two separate integrals. However, in the calculations that follow later, the two residues turn out to be additive inverses of each other and hence do not contribute anything.

The following is an analogue of Ramanujan's transformation formula (Theorem \ref{entry1}) for odd characters.
\begin{corollary}\label{1splcaseso}
For an odd character $\chi$, define $Q(\alpha,\chi)$ by
\begin{equation*}
Q(\alpha,\chi)=\sqrt{\alpha}\hspace{1mm} \textup{Im}\left(G(\chi)\sum_{n=1}^{\infty}\sum_{k=1}^{\infty}\frac{\overline{\chi}(n)\overline{\chi}(k)}{k+n\alpha}\right)=-\sqrt{\alpha}\hspace{1mm} \textup{Im}\left(G(\chi)\sum_{n=1}^{\infty}\overline{\chi}(n)\psi\left(n\alpha,\overline{\chi}\right)\right),
\end{equation*}
where $\psi\left(a,\chi\right)$ is defined in \textup{(\ref{charpsi})}. Then we have
\begin{align}\label{splcaseso}
&Q(\alpha,\chi)=Q(\beta,\overline{\chi})=Q(\alpha,\overline{\chi})=Q(\beta,\chi)\nonumber\\
&=\frac{1}{4\pi^{1/2}q}\int_{0}^{\infty}\Gamma\bigg(\frac{1+it}{4}\bigg)\Gamma\bigg(\frac{1-it}{4}\bigg)\Xi\left(\frac{t}{2},\overline{\chi}\right)\Xi\left(\frac{t}{2},\chi\right)\cos\left(\frac{1}{2}t\log\alpha\right)dt.
\end{align}
\end{corollary}
\begin{proof}
Using Lemma \ref{asycharpsi}, we find that the double series in the definition of $Q(\alpha,\chi)$ converges. Let $z\to 0$ in Theorem \ref{ramgenoc}. Multiplying both sides by $-q$ and using (\ref{sgs}) and (\ref{eo}), we observe that
\begin{align}\label{splcaseo}
&\sqrt{\alpha}\left(G(\overline{\chi})\sum_{n=1}^{\infty}\sum_{k=1}^{\infty}\frac{\chi(n)\chi(k)}{k+n\alpha}+G(\chi)\sum_{n=1}^{\infty}\sum_{k=1}^{\infty}\frac{\overline{\chi}(n)\overline{\chi}(k)}{k+n\alpha}\right)\nonumber\\
&=\sqrt{\beta}\left(G(\chi)\sum_{n=1}^{\infty}\sum_{k=1}^{\infty}\frac{\overline{\chi}(n)\overline{\chi}(k)}{k+n\beta}+G(\overline{\chi})\sum_{n=1}^{\infty}\sum_{k=1}^{\infty}\frac{\chi(n)\chi(k)}{k+n\beta}\right)\nonumber\\
&=\frac{i}{2\pi^{1/2}q}\int_{0}^{\infty}\Gamma\bigg(\frac{1+it}{4}\bigg)\Gamma\bigg(\frac{1-it}{4}\bigg)\Xi\left(\frac{t}{2},\overline{\chi}\right)\Xi\left(\frac{t}{2},\chi\right)\cos\left(\frac{1}{2}t\log\alpha\right)\, dt.
\end{align}
Now using (\ref{eo}) for odd characters to simplify (\ref{splcaseo}), we see that
\begin{align*}\label{innovo}
&2i\sqrt{\alpha}\hspace{1mm}\textup{Im} \left(G(\overline{\chi})\sum_{n=1}^{\infty}\sum_{k=1}^{\infty}\frac{\chi(n)\chi(k)}{k+n\alpha}\right)=2i\sqrt{\alpha}\hspace{1mm} \textup{Im} \left(G(\chi)\sum_{n=1}^{\infty}\sum_{k=1}^{\infty}\frac{\overline{\chi}(n)\overline{\chi}(k)}{k+n\alpha}\right)\nonumber\\
&=2i\sqrt{\beta}\hspace{1mm}\textup{Im}\left(G(\chi)\sum_{n=1}^{\infty}\sum_{k=1}^{\infty}\frac{\overline{\chi}(n)\overline{\chi}(k)}{k+n\beta}\right)=2i\sqrt{\beta}\hspace{1mm} \textup{Im} \left(G(\overline{\chi})\sum_{n=1}^{\infty}\sum_{k=1}^{\infty}\frac{\chi(n)\chi(k)}{k+n\beta}\right)\nonumber\\
&=\frac{i}{2\pi^{1/2}q}\int_{0}^{\infty}\Gamma\bigg(\frac{1+it}{4}\bigg)\Gamma\bigg(\frac{1-it}{4}\bigg)\Xi\left(\frac{t}{2},\overline{\chi}\right)\Xi\left(\frac{t}{2},\chi\right)\cos\left(\frac{1}{2}t\log\alpha\right)\, dt.
\end{align*}
This implies (\ref{splcaseso}).
\end{proof}
If we now start with the integral in Corollary \ref{1splcaseso}, evaluate it using (\ref{gfica1}) with $z=0$ and make use of Lemma \ref{invmelimp} when $\chi$ is odd, we obtain the same result as in Corollary \ref{1splcaseso}, except that the function $Q(\alpha,\chi)$ is replaced by $-iF(\alpha,\chi)$, where $F(\alpha,\chi)$ is defined in (\ref{ta}). It is then trivial to see that $F(\alpha,\chi)=iQ(\alpha,\chi)$.

We separately record the following theorem resulting from the discussion on the previous line and the one succeeding Corollary \ref{1splcases}.
\begin{theorem}\label{rim}
The sum $F(\alpha,\chi)$ defined in \textup{(\ref{ta})} is real if $\chi$ is even and purely imaginary if $\chi$ is odd.
\end{theorem}

\section{Character analogues of the Ramanujan-Hardy-Littlewood conjecture}
In this section, we prove Theorems \ref{rhlao} and \ref{rhlae}. We require Lemma 3.1 from \cite{abyz} which states that if $\chi$ is a primitive character of conductor $N$ and $k\geq 2$ is an integer such that $\chi(-1)=(-1)^{k}$, then
\begin{equation}\label{fabyz}
\frac{(k-2)!N^{k-2}G(\chi)}{2^{k-1}\pi^{k-2}i^{k-2}}L(k-1,\overline{\chi})=L'(2-k,\chi).
\end{equation} 
\begin{proof}[Theorem \textup{\ref{rhlao}}][]
From \cite{landau}, we have for Re $s>1$,
\begin{equation}\label{l2}
\sum_{n=1}^{\infty}\frac{\chi(n)\mu(n)}{n^s}=\frac{1}{L(s,\chi)}.
\end{equation}
Also, since for $-1<c=$ Re $s<0$,
\begin{equation}\label{gme}
(1-e^{-x})=-\frac{1}{2\pi i}\int_{c-i\infty}^{c+i\infty}\Gamma(s)x^{-s}\, ds,
\end{equation}
replacing $s$ by $s+1$, we find that for $-2<c<-1$,
\begin{equation}\label{gm2}
(1-e^{-x})=-\frac{1}{2\pi i}\int_{c-i\infty}^{c+i\infty}\Gamma(s+1)x^{-s-1}\, ds.
\end{equation}
Using (\ref{l2}) and (\ref{gm2}), we observe that
{\allowdisplaybreaks\begin{align}\label{icalcu}
\sum_{n=1}^{\infty}\frac{\chi(n)\mu(n)}{n^2}e^{-\frac{\pi\alpha^2}{n^2q}}&=\frac{1}{L(2,\chi)}-\sum_{n=1}^{\infty}\frac{\chi(n)\mu(n)}{n^2}(1-e^{-\frac{\pi\alpha^2}{n^2q}})\nonumber\\
&=\frac{1}{L(2,\chi)}+\frac{q}{2\pi^2 i\alpha^2}\int_{c-i\infty}^{c+i\infty}\sum_{n=1}^{\infty}\frac{\chi(n)\mu(n)}{n^{-2s}}\Gamma(s+1)\left(\frac{\pi\alpha^{2}}{q}\right)^{-s}\, ds\nonumber\\
&=\frac{1}{L(2,\chi)}+\frac{q}{2\pi^2 i\alpha^2}\int_{c-i\infty}^{c+i\infty}\frac{\Gamma(s+1)}{L(-2s,\chi)}\left(\frac{\pi\alpha^{2}}{q}\right)^{-s}\, ds,
\end{align}}%
where in the second step above, we interchanged the order of summation and integration because of absolute convergence. For $\chi$ odd, (\ref{funl}) can be put in the form
\begin{equation*}
\left(\frac{\pi}{q}\right)^{-(2-s)/2}\Gamma\left(\frac{2-s}{2}\right)L(1-s,\overline{\chi})=\frac{iq^{1/2}}{G(\chi)}\left(\frac{\pi}{q}\right)^{-(s+1)/2}\Gamma\left(\frac{s+1}{2}\right)L(s,\chi).
\end{equation*}
Hence,
\begin{equation}\label{oddfu}
\frac{\Gamma(s+1)}{L(-2s,\chi)}=\frac{G(\overline{\chi})}{iq^{1/2}}\left(\frac{\pi}{q}\right)^{2s+\tf{1}{2}}\frac{\Gamma\left(\frac{1}{2}-s\right)}{L(2s+1,\overline{\chi})}.
\end{equation}
Substituting (\ref{oddfu}) in (\ref{icalcu}), we observe that
\begin{equation}\label{icalcu1}
\sum_{n=1}^{\infty}\frac{\chi(n)\mu(n)}{n^2}e^{-\frac{\pi\alpha^2}{n^2q}}=\frac{1}{L(2,\chi)}-\frac{G(\overline{\chi})}{2\pi^{3/2}\alpha^2}\int_{c-i\infty}^{c+i\infty}\frac{\Gamma\left(\frac{1}{2}-s\right)}{L(2s+1,\overline{\chi})}\left(\frac{q\alpha^{2}}{\pi}\right)^{-s}\, ds.
\end{equation}
We wish to shift the line of integration from Re $s=c$, $-2<c<-1$, to Re $s=\lambda$, where $\frac{1}{2}<\lambda<\frac{3}{2}$. Consider a positively oriented rectangular contour formed by $[c-iT, \lambda-iT], [\lambda-iT, \lambda+iT], [\lambda+iT,c+iT]$ and $[c+iT,c-iT]$, where $T$ is any positive real number. Let $\rho=\delta+i\gamma$ denote a non-trivial zero of $L(s,\overline{\chi})$. Let $T\to\infty$ through values such that $|T-\gamma|>\exp\left(-A_{1}\gamma/\log\gamma\right)$ for every ordinate $\gamma$ of a zero of $L(s,\overline{\chi})$. It is known \cite[p.~102]{dav} that for $t$ not coinciding with the ordinate $\gamma$ of a zero, and $-1\leq\sigma\leq 2$, 
\begin{equation*}
\frac{L^{'}(s,\overline{\chi})}{L(s,\overline{\chi})}=\sum_{|t-\gamma|\leq 1}\frac{1}{s-\rho}+O\left(\log\left(q(|t|+2)\right)\right),
\end{equation*}
from which we can conclude that
\begin{equation}\label{nb}
\log L(s,\overline{\chi})=\sum_{|t-\gamma|\leq 1}\log(s-\rho)+O\left(\log\left(q(|t|+2)\right)\right).
\end{equation}
Taking real parts in (\ref{nb}) gives
{\allowdisplaybreaks\begin{align}\label{nb1}
\log |L(s,\overline{\chi})|&=\sum_{|t-\gamma|\leq 1}\log|s-\rho|+O\left(\log\left(q(|t|+2)\right)\right)\nonumber\\
&\geq \sum_{|t-\gamma|\leq 1}\log|t-\gamma|+O\left(\log\left(q(|t|+2)\right)\right).\nonumber\\
\end{align}}
Hence from (\ref{nb1}), we have
\begin{align}\label{nb12}
\log |L(\sigma+iT,\overline{\chi})|&\geq -\sum_{|T-\gamma|\leq 1}A_{1}\gamma/\log\gamma+O\left(\log\left(q(|T|+2)\right)\right)\nonumber\\
&>-A_{2}T,
\end{align}
where $A_{2}<\pi/4$ if $A_{1}$ is small enough and $T>T_{0}$ for some fixed $T_{0}$. From (\ref{nb12}), we see that
\begin{equation}\label{nb2}
\left|\frac{1}{L(2s+1,\overline{\chi})}\right|<e^{A_{3}T},
\end{equation}
where $A_{3}<\pi/2$. Using (\ref{strivert}) and (\ref{nb2}), we observe that as $T\to\infty$ through the above values, the integrals along the horizontal segments tend to zero. Now let $\frac{\rho-1}{2}:=\delta+i\gamma$ denote a non-trivial zero of $L(2s+1,\overline{\chi})$. Let $R_{f}(a)$ denote the residue at $a$ of the function $f(s):=\displaystyle\frac{\Gamma\left(\frac{1}{2}-s\right)}{L(2s+1,\overline{\chi})}\left(\frac{q\alpha^{2}}{\pi}\right)^{-s}$. The non-trivial zeros of $L(2s+1,\overline{\chi})$ lie in the critical strip $-\frac{1}{2}<$ Re $s<0$, whereas the trivial zeros are at $-1, -2, -3,\cdots$. Also, $\Gamma\left(\frac{1}{2}-s\right)$ has poles at $\frac{1}{2}, \frac{3}{2}, \frac{5}{2},\cdots$. Then the residue theorem yields
{\allowdisplaybreaks\begin{align}\label{resmra}
\int_{c-i\infty}^{c+i\infty}\frac{\Gamma\left(\frac{1}{2}-s\right)}{L(2s+1,\overline{\chi})}\left(\frac{q\alpha^{2}}{\pi}\right)^{-s}\, ds&=\int_{\lambda-i\infty}^{\lambda+i\infty}\frac{\Gamma\left(\frac{1}{2}-s\right)}{L(2s+1,\overline{\chi})}\left(\frac{q\alpha^{2}}{\pi}\right)^{-s}\, ds\nonumber\\
&\quad-2\pi i\left(R_{f}(-1)+\sum_{\rho}R_{f}\left(\frac{\rho-1}{2}\right)+R_{f}\left(\frac{1}{2}\right)\right),
\end{align}}
where
\begin{align}\label{reszmra}
R_{f}(-1)&=\lim_{s\to -1}(s+1)\frac{\Gamma\left(\frac{1}{2}-s\right)}{L(2s+1,\overline{\chi})}\left(\frac{q\alpha^{2}}{\pi}\right)^{-s}=\frac{\alpha^{2}q}{4\sqrt{\pi}L'(-1,\overline{\chi})},\nonumber\\
R_{f}\left(\frac{\rho-1}{2}\right)&=\lim_{s\to\frac{\rho-1}{2}}\left(s-\frac{\rho-1}{2}\right)\frac{\Gamma\left(\frac{1}{2}-s\right)}{L(2s+1,\overline{\chi})}\left(\frac{q\alpha^{2}}{\pi}\right)^{-s}=\frac{\Gamma\left(\frac{2-\rho}{2}\right)}{2L'(\rho,\overline{\chi})}\left(\frac{\pi}{q\alpha^{2}}\right)^{\frac{\rho-1}{2}},\nonumber\\
R_{f}(1/2)&=-\frac{\sqrt{\pi}}{\alpha\sqrt{q}L(2,\overline{\chi})}.
\end{align}
Of course, here we have assumed that the non-trivial zeros of $L(2s+1,\overline{\chi})$ are all simple and that $\sum_{\rho}R_{f}\left(\frac{\rho-1}{2}\right)$ converges, since the afore-mentioned discussion regarding the integrals along the horizontal segments tending to zero as $T\to\infty$ through the chosen sequence does not imply the convergence of $\sum_{\rho}R_{f}\left(\frac{\rho-1}{2}\right)$ in the ordinary sense. It only means that the series converges only when we bracket the terms in such a way that the two terms for which
\begin{equation*}
|\gamma-\gamma'|<\exp\left(-A_{1}\gamma/\log 2\gamma\right)+\exp\left(-A_{1}\gamma'/\log 2\gamma'\right)
\end{equation*}
are included in the same bracket.
Using (\ref{l2}) and interchanging the order of summation and integration because of absolute convergence, we obtain
\begin{align}\label{lamint}
\int_{\lambda-i\infty}^{\lambda+i\infty}\frac{\Gamma\left(\frac{1}{2}-s\right)}{L(2s+1,\overline{\chi})}\left(\frac{q\alpha^{2}}{\pi}\right)^{-s}\, ds&=\sum_{n=1}^{\infty}\frac{\overline{\chi}(n)\mu(n)}{n}\int_{\lambda-i\infty}^{\lambda+i\infty}\Gamma\left(\frac{1}{2}-s\right)\left(\frac{q\alpha^{2}n^2}{\pi}\right)^{-s}\, ds\nonumber\\
&=\frac{\sqrt{\pi}}{\alpha\sqrt{q}}\sum_{n=1}^{\infty}\frac{\overline{\chi}(n)\mu(n)}{n^2}\int_{d-i\infty}^{d+i\infty}\Gamma(s)\left(\frac{\pi}{\alpha^{2}n^2q}\right)^{-s}\, ds,
\end{align}
where in the penultimate line, we have made the change of variable $s\to\frac{1}{2}-s$ so that $-1<d<0$. Thus, (\ref{resmra}), (\ref{reszmra}), (\ref{lamint}) and (\ref{gme}) imply
\begin{align}\label{icalcu2}
&\int_{c-i\infty}^{c+i\infty}\frac{\Gamma\left(\frac{1}{2}-s\right)}{L(2s+1,\overline{\chi})}\left(\frac{q\alpha^{2}}{\pi}\right)^{-s}\, ds\nonumber\\
&=-\frac{2\pi^{3/2} i}{\alpha\sqrt{q}}\sum_{n=1}^{\infty}\frac{\overline{\chi}(n)\mu(n)}{n^2}\left(1-e^{-\frac{\pi}{\alpha^2n^2q}}\right)\nonumber\\
&\quad-2\pi i\left(\frac{\alpha^{2}q}{4\sqrt{\pi}L'(-1,\overline{\chi})}+\sum_{\rho}\frac{\Gamma\left(\frac{2-\rho}{2}\right)}{2L'(\rho,\overline{\chi})}\left(\frac{\pi}{q\alpha^{2}}\right)^{\frac{\rho-1}{2}}-\frac{\sqrt{\pi}}{\alpha\sqrt{q}L(2,\overline{\chi})}\right).
\end{align}
From (\ref{icalcu1}), (\ref{icalcu2}) and the facts that $\alpha\beta=1$ and $\sqrt{G(\chi)G(\overline{\chi})}=i\sqrt{q}$, we find that
{\allowdisplaybreaks\begin{align}\label{icalcuq}
&\alpha\sqrt{\alpha}\sqrt{G(\chi)}\sum_{n=1}^{\infty}\frac{\chi(n)\mu(n)}{n^2}e^{-\frac{\pi\alpha^2}{n^2q}}\nonumber\\
&=\frac{\alpha\sqrt{\alpha}\sqrt{G(\chi)}}{L(2,\chi)}-\frac{\beta\sqrt{\beta}\sqrt{G(\overline{\chi})}}{L(2,\overline{\chi})}+\beta\sqrt{\beta}\sqrt{G(\overline{\chi})}\sum_{n=1}^{\infty}\frac{\overline{\chi}(n)\mu(n)}{n^2}e^{-\frac{\pi\beta^2}{n^2q}}\nonumber\\
&\quad-\frac{\alpha\sqrt{\alpha}q^{3/2}\sqrt{G(\overline{\chi})}}{4\pi L'(-1,\overline{\chi})}-\frac{q\sqrt{G(\overline{\chi})}}{2\pi\sqrt{\beta}}\sum_{\rho}\frac{\Gamma\left(\frac{2-\rho}{2}\right)}{L'(\rho,\overline{\chi})}\left(\frac{\pi}{q}\right)^{\rho/2}\beta^{\rho}+\frac{\beta\sqrt{\beta}\sqrt{G(\overline{\chi})}}{L(2,\overline{\chi})}.
\end{align}}
Applying (\ref{fabyz}) with $N=q$ and $k=3$ and replacing $\chi$ by $\overline{\chi}$ gives
\begin{equation}\label{derl}
\frac{1}{L'(-1,\overline{\chi})}=\frac{4\pi i}{qG(\overline{\chi})L(2,\chi)}.
\end{equation}
Thus (\ref{icalcuq}) and (\ref{derl}) yield
\begin{align}\label{intrhl}
&\alpha\sqrt{\alpha}\sqrt{G(\chi)}\sum_{n=1}^{\infty}\frac{\chi(n)\mu(n)}{n^2}e^{-\frac{\pi\alpha^2}{n^2q}}-\beta\sqrt{\beta}\sqrt{G(\overline{\chi})}\sum_{n=1}^{\infty}\frac{\overline{\chi}(n)\mu(n)}{n^2}e^{-\frac{\pi\beta^2}{n^2q}}\nonumber\\
&=-\frac{q\sqrt{G(\overline{\chi})}}{2\pi\sqrt{\beta}}\sum_{\rho}\frac{\Gamma\left(\frac{2-\rho}{2}\right)}{L'(\rho,\overline{\chi})}\left(\frac{\pi}{q}\right)^{\rho/2}\beta^{\rho}.
\end{align}
Switching the roles of $\alpha$ and $\beta$ and those of $\chi$ and $\overline{\chi}$ gives
\begin{equation}\label{intrhl1}
\frac{q\sqrt{G(\chi)}}{2\pi\sqrt{\alpha}}\sum_{\rho}\frac{\Gamma\left(\frac{2-\rho}{2}\right)}{L'(\rho,\chi)}\left(\frac{\pi}{q}\right)^{\rho/2}\alpha^{\rho}+\frac{q\sqrt{G(\overline{\chi})}}{2\pi\sqrt{\beta}}\sum_{\rho}\frac{\Gamma\left(\frac{2-\rho}{2}\right)}{L'(\rho, \overline{\chi})}\left(\frac{\pi}{q}\right)^{\rho/2}\beta^{\rho}=0.                 
\end{equation}
Finally (\ref{intrhl}) and (\ref{intrhl1}) give (\ref{mrao}) upon simplification.
\end{proof}
\textbf{Remark.} The approach used above for proving that the integrals along the horizontal segments tend to zero as $T\to\infty$ through the chosen sequence is adapted from \cite[p.~219]{titch}.

To prove Theorem \ref{rhlae}, we require the following lemma.
\begin{lemma}\label{kl}
\begin{equation*}
\sum_{n=1}^{\infty}\frac{\chi(n)\mu(n)}{n}=\frac{1}{L(1,\chi)}.
\end{equation*}
\end{lemma}
\begin{proof}
Dividing $n$ into its residue classes mod $q$ by letting $n=qr+b$, $0\leq r<\infty$, $0\leq b\leq q-1$, we find that since $\chi$ has period $q$,
\begin{equation}\label{kl1}
\sum_{n=1}^{\infty}\frac{\chi(n)\mu(n)}{n}=\sum_{r=0}^{\infty}\sum_{b=0}^{q-1}\frac{\chi(b)\mu(qr+b)}{qr+b}=\sum_{b=0}^{q-1}\chi(b)\sum_{r=0}^{\infty}\frac{\mu(qr+b)}{qr+b}.
\end{equation}
The series $\sum_{r=0}^{\infty}\mu(qr+b)/(qr+b)$ was first studied by J.C.~Kluyver \cite{kluyver} and its convergence was proved by E.~Landau \cite{landau}. In fact, Landau gave an explicit representation for this series in terms of a finite sum consisting of $L$-functions. Thus (\ref{kl1}) implies convergence of $\sum_{n=1}^{\infty}\chi(n)\mu(n)/n$. Then using (\ref{l2}) and an analogue of Abel's theorem for power series, we see that 
\begin{equation*}
\sum_{n=1}^{\infty}\frac{\chi(n)\mu(n)}{n}=\lim_{s\to 1}\sum_{n=1}^{\infty}\frac{\chi(n)\mu(n)}{n^s}=\lim_{s\to 1}\frac{1}{L(s,\chi)}=\frac{1}{L(1,\chi)}.
\end{equation*}
\end{proof}
\begin{proof}[Theorem \textup{\ref{rhlae}}][]
The proof is very similar to that of Theorem \ref{rhlao} and hence we omit the details. 
However we note that Lemma \ref{kl}, (\ref{funl}) in the form \cite[p.~69]{dav}
\begin{equation*}
\pi^{-(1-s)/2}q^{(1-s)/2}\Gamma\left(\frac{1-s}{2}\right)L(1-s,\overline{\chi})=\frac{q^{1/2}}{G(\chi)}\pi^{-s/2}q^{s/2}\Gamma\left(\frac{s}{2}\right)L(s,\chi)
\end{equation*}
and (\ref{fabyz}) with $N=q$ and $k=2$ are used in the proof. 
\end{proof}
\section{Open problems}
The following are some open problems with which we will conclude.

1. We have indirectly given the proof of the fact that function $F(\alpha,\chi)$ defined in (\ref{ta}) is real (purely imaginary) respectively when $\chi$ is even (odd). Prove this directly, i.e., without using Corollaries \ref{1splcases} and \ref{1splcaseso} and the integrals in those corollaries.

2. Since (\ref{mr}) is of the form $F(\alpha)=F(\beta)$, where $\alpha\beta=1$, it is natural to ask if there exists an integral representation involving the Riemann $\Xi$-function equal to the two expressions in (\ref{mr}). Finding an integral representation for either side of (\ref{mr}) may throw light on the convergence of $\sum_{\rho}\left(\Gamma\left((1-\rho)/2\right)a^{\rho}\right)/\zeta^{'}(\rho)$ provided, of course, that the integral converges in the first place.
It should be remarked here that Hardy and Littlewood \cite[p.~161]{hl} have shown that the relation 
\begin{equation}\label{reisz}
P(y)=O\left(y^{-\frac{1}{4}+\delta}\right),
\end{equation}
where $P(y)=\sum_{n=1}^{\infty}(-y)^{n}/(n!\zeta(2n+1))$, can be derived from (\ref{mr}) provided we assume the Riemann hypothesis and the absolute convergence of $\sum_{\rho}\Gamma\left((1-\rho)/2\right)/\zeta^{'}(\rho)$. They have further shown that (\ref{reisz}) is a necessary and sufficient condition for the Riemann hypothesis to be true.

Similarly, it is natural to ask if the expressions in (\ref{mrao}) and (\ref{mrae}) have integral representations involving $\Xi\left(\frac{t}{2},\chi\right)$.\\

3. Does there exist a generalization of Theorem \ref{rhl} admitting representations of the form $F(z,\alpha)=F(z,\beta)$? Similarly, do there exist generalizations of Theorems \ref{rhlao} and \ref{rhlae} admitting representations of the form $F(z, \alpha,\chi)=F(-z, \beta,\overline{\chi})=F(-z,\alpha,\overline{\chi})=F(z,\beta,\chi)$?

\begin{center}
\textbf{Acknowledgements}
\end{center}
The author wishes to express his sincere thanks to Professor Bruce C.~Berndt for his constant support, careful reading of this manuscript and for several helpful comments, and to Professors Alexandru Zaharescu, Kevin Ford and Adolf Hildebrand for their help. He would also like to thank Becky S.~Burner, Timothy W.~Cole and Margaret A.~Lewis from the Mathematics library at University of Illinois at Urbana-Champaign, and Anton Lukyanenko for their great help in locating the Russian references and citations. Last but not the least, he would like to thank Reetesh Ranjan for numerically verifying some of the identities by programming in Fortran.

\end{document}